\documentclass[10pt]{article}
 \usepackage[T1]{fontenc}
\usepackage{natbib}
\usepackage{amsmath,amsthm, amssymb}
 \usepackage{graphicx}
 \usepackage{dsfont}
 \usepackage{color}
 \usepackage{float}

\renewcommand{\L}{\mathbb{L}}
\newcommand{\vvvert}{\rvert\hspace{-0.12em}\rvert\hspace{-0.12em}\rvert}

\newcommand{\E}{\ensuremath{\mathbb{E}}}
\renewcommand{\P}{\ensuremath{\mathbb{P}}}
\newcommand{\R}{\ensuremath{\mathbb{R}}}
\newcommand{\fx}{\mathrm{f}_X}
\newcommand{\hfx}{\hat{\mathrm{f}}_X}
\newcommand{\e}{\ensuremath{\varepsilon}}

\newcommand{\var}{\mbox{var}}

\newcommand{\red}[1]{{\color{red}{\bf #1}}}
\newcommand{\pp}{q}
\newcommand{\argmin}{\mathop{\mathrm{arg\,min}}}
\def\1{\mathds{1}}
\def\Sp{\mathrm{Sp}}

\numberwithin{equation}{section}

\newtheorem{definition}{Definition}

\newtheorem{proposition}{Proposition}
 \newtheorem{lemma}{Lemma}
\newtheorem{theorem}{Theorem}
\newtheorem{cor}{Corollary}

\begin{document}

\markboth{\hfill{\footnotesize\rm KARINE BERTIN AND CLAIRE LACOUR AND VINCENT RIVOIRARD
}\hfill}
{\hfill {\footnotesize\rm Adaptive estimation of conditional density function} \hfill}
\renewcommand{\thefootnote}{}
$\ $\par

\centerline{\large\bf Adaptive pointwise estimation of conditional density function}
\vspace{.4cm}
\centerline{Karine Bertin, Claire Lacour and Vincent Rivoirard}
\vspace{.4cm}
\centerline{\it  Universidad de Valpara\'{i}so, Universit\'e Paris-Sud, Universit\'e Paris Dauphine
}
\vspace{.55cm}


\begin{quotation}
\noindent {\it Abstract:}
In this paper we consider the problem of estimating $f$, the conditional density  of $Y$ given $X$, by using an independent sample distributed as $(X,Y)$ in the multivariate setting. We consider the estimation of $f(x,.)$ where $x$ is a fixed point. We define two different procedures of estimation, the first one using kernel rules, the second one inspired from projection methods. Both adapted estimators are tuned by using the Goldenshluger and Lepski  methodology. After deriving lower bounds, we show that these procedures satisfy oracle inequalities and are optimal from the minimax point of view on anisotropic H\"{o}lder balls. Furthermore, our results allow us to measure precisely the influence of $\fx(x)$ on rates of convergence, where $\fx$ is the density of $X$.  Finally, some simulations illustrate the good behavior of our tuned estimates in practice.
\par

\vspace{9pt}
\par
\end{quotation}\par

\section{Introduction}
\subsection{Motivation}\label{motivation}
In this paper, we consider the problem of conditional density estimation. For this purpose, we assume we are given an i.i.d. sample $(X_i,Y_i)$ of couples of random vectors (for any $i$, $X_i\in\R^{d_1}$ and $Y_i\in\R^{d_2}$, with $d_1\geq 1$ and $d_2\geq 1$) with common
probability density function $f_{X,Y}$ and marginal densities $f_{Y}$ and  $\fx$: for any  $y\in\R^{d_2}$  and any $x\in\R^{d_1}$,
$$f_{Y}(y)=\int_{\R^{d_1}}f_{X,Y}(u,y)du,\quad \fx(x)=\int_{\R^{d_2}}f_{X,Y}(x,v)dv.$$
The conditional density function of $Y_i$ given $X_i = x$ is defined by
$$f(x,y)=\frac{f_{X,Y}(x,y)}{\fx(x)}$$
for all $y\in\mathbb{R}^{d_2}$ and
$x\in\mathbb{R}^{d_1}$ such that $\fx (x) > 0.$
Our goal is to estimate $f$ using the
observations $(X_i,Y_i)$.  The conditional density is much more informative than the simple regression function and then its estimation has many practical applications: in Actuaries (\cite{efromovich10}), Medicine (\cite{takeuchi}), Economy  (\cite{HallRacineLi}), Meteorology (\cite{JeonTaylor}) among others. In particular, due to recent advances in ABC methods, the problem of conditional density estimation in the multivariate setting is of main interest.

Indeed, the ABC methodology, where ABC stands for {\em approximate Bayesian computation},
offers a resolution of untractable-yet-simulable models, that is models for which it is impossible
to calculate the likelihood.
The standard ABC procedure is very intuitive and consists in
\begin{itemize}
\item simulating a lot of parameters values using the prior distribution and, for each parameter
value, a corresponding dataset,
\item comparing this simulated dataset to the observed one;
\item finally, keeping the parameter values for which distance between the simulated dataset and
the observed one is smaller than a tolerance level.
\end{itemize}
That is a crude nonparametric approximation of the target posterior distribution (the conditional
distribution of the parameters given the observation).
Even if some nonparametric perspectives have been considered (see
\cite{blum:2010} or \cite{biauetal2012}), we easily imagine that, using the simulated couples (parameters and datasets), a good nonparametric estimation of
the posterior distribution can be a credible alternative to the ABC method.
Such a procedure has to consider that the conditional density has to be estimated only for the observed
value in the conditioning.

All previous points clearly motivate our work and in the sequel, we aim at providing an estimate with the following 4 requirements:
\begin{enumerate}
\item The estimate has to be fully data-driven and implementable in a reasonable computational time.
\item  The parameters of the method have to adapt to the function $f$ in the neighborhood of $x$. Tuning the hyperparameters of the estimate has to be an easy task. 
\item The estimate should be optimal from the theoretical point of view in an asymptotic setting but also in a non-asymptotic one.
\item Estimating $f$ in neighborhoods of points $x$ where $\fx(x)$ is equal or close to 0  is of course a difficult task and a loss is unavoidable. Studying this loss  and providing estimates that are optimal with respect to this problem are the fourth motivation of this paper.
\end{enumerate}

To address the problem of  conditional density estimation, the first idea of statisticians was to estimate $f$ by the ratio of a kernel estimator of the joint density $f_{X,Y}$ and a kernel estimator of $\fx$: see \cite{Rosenblatt69}, \cite{ChenLintonRobinson}, or also \cite{HBG96}, \cite{GZ} for refinements of this method.
A important work in this line is the one of  \cite{FYT96} who extend the Rosenblatt estimator by a local polynomial method (see also \cite{HY02}).
The estimators introduced in the ABC literature are also of this kind: a linear (or quadratic) adjustment is realized on the data before applying the classic quotient estimator (\cite{beaumont2002} , \cite{blum:2010}).
Other directions are investigated by
\cite{BouazizLopez} who use a single-index model, or \cite{GK} who partition the space and obtain a piecewise constant estimate.
All these papers have in common to involve a ratio between two density estimates, though we can mention \cite{stone94} for a spline tensor based maximum likelihood estimator.
An original approach which rather involves a product is the copula one of \cite{faugeras}. But his method depends on a bandwidth, that remains to select from the data. In particular, for all of these methods, the second requirement is not satisfied.

The practical choice of the bandwidth and cross-validation methods are studied in \cite{BH01} and \cite{FanYim04}.
However, no theoretical result is associated to this study. The first adaptive results can be found in \cite{clem} for the estimation of the transition density of a Markov chain, which is a very similar problem to the one of conditional density estimation (set $Y_{i}=X_{i+1}$). He uses thresholding of wavelet estimator. Afterwards, using different methods,  the works of \cite{brunelcomtelacour07} or \cite{efromovich07} yield oracle inequalities and minimax rates of convergence for anisotropic conditional densities. The case of inhomogenous regularities is studied in \cite{akakpolacour} or  \cite{Sart13} in the case of Markov chains.
Still for global adaptive approach, we can cite \cite{Chagny} who applies the Goldenshluger-Lepski methodology to warped bases and  \cite{CohenLePennec} who use a model selection approach with Kullback risk. All the previous authors use a global risk and either consider integration with respect to $\fx(x)dx$ or assume that $\fx$ is bounded from below by a constant (as it is done in regression estimation). We are interested in precisely studying this assumption to show that it is unavoidable in  some sense.

\subsection{Our strategy and our contributions}\label{sec:strategy}
Our strategy to estimate $f$ is based on the Goldenshluger and Lepski methodology proposed in the seminal papers \cite{GL11, GL6} in the case of density estimation and extended to the white noise and regression models in \cite{GL5}. This strategy detailed in Section~\ref{explicationGLM} allows us to derive two procedures: kernel and projection rules. If they seem different, they are based on similar ideas and they lead to quite similar theoretical results. Our method automatically selects a regularization parameter, and in particular a bandwidth for kernel rules. Note that the tolerance level in ABC methods can be reinterpreted as a regularization parameter.

Unlike most of previous works of the literature, we shall not use a global risk and we will evaluate the quality of an estimator $\hat f$ at a fixed point $x\in\mathbb{R}$ and in the $\L_2$-norm with respect to the variable $y$. In other words, we will use the risk
\begin{equation}\label{riskx2}
R_x(\hat f,q)=\left(\E\left[\|\hat f-f\|_{x,2}^\pp \right]\right)^{\frac{1}{q}},
\end{equation}
where for any function $g$, 
\begin{equation}\label{normx2}
\|g\|_{x,2}=\left(\int_{\mathbb{R}^{d_2}}g^2(x,y)dy\right)^{1/2}.
\end{equation}
The  previously mentioned motivating applications show that the tuning parameter has to depend on $x$, which is not the case of other cited-above adaptive methods. As shown later, combined with the Goldenshluger and Lepski methodology, considering this risk allows us to derive estimates satisfying this property. Furthermore, for a given $x$, $y\mapsto f(x,y)$ is a density, so it is natural for us to study the estimation pointwisely in $x$.

From the theoretical point of view, we establish non asymptotic meaningful oracle inequalities and rates of convergence on anisotropic H\"{o}lder balls $\mathcal{H}_d(\boldsymbol{\alpha},\bold L)$. More precisely, in Proposition~\ref{borneinforacle} and Theorem~\ref{t2},  we establish lower bounds in oracle and minimax settings. Then, upper bounds of the risk for our adaptive kernel procedure are established (see Theorems~\ref{theo1}, \ref{theo2} and \ref{propvitesse}).
 If the density $\fx$ is smooth enough, Corollary~\ref{cor} shows that upper and lower bounds match up to constants in the asymptotic setting. Then, there is a natural question: is this assumption on the smoothness  of $\fx$ mandatory? We prove that the answer is no by establishing the upper bound of the risk for our adaptive projection estimate (see Theorems~\ref{oracleproba} and \ref{vitesse}). In particular, the latter achieves a polynomial rate of convergence on anisotropic H\"{o}lder balls with rate exponent $\bar\alpha/(2\bar\alpha+1)$,
where $\bar \alpha$ is the classical anisotropic smoothness
index.
To our knowledge, this rate exponent is new in the conditional density estimation setting for the pointwise risk in $x$. Our result also explicits the dependence of the rate with respect to $\bold{L}$ on the one hand and to $\fx(x)$ on the other hand, which is not classical. Indeed,  as previously recalled, estimation is harder when $\fx(x)$ is small and this is the reason why most of the papers assume that $\fx$ is bounded from below by a constant. For kernel rules, our study is sharp enough to measure precisely the influence of $\fx(x)$ on  the performance of our procedure. Under some conditions and if the sample size is $n$, we show that the order of magnitude of minimax rates (that are achieved by our procedure), is $(n\fx(x))^{\bar\alpha/(2\bar\alpha+1)}$. We conclude that our setting is equivalent to the setting where $\fx$ is locally bounded from $0$ by 1 but we observe $n\fx(x)$ observations instead of $n$.

Finally, we study our procedures from a practical point of view. We aim at completing theoretical results by studying tuning issues. More precisely, our procedures are data driven and tuning parameters depend on $x$ and on an hyperparameter $\eta$, a constant that has to be tuned. We lead a precise study that shows how to choose $\eta$ in practice. We also show that reconstructions for various examples and various values of $n$ are satisfying.  All these results show that our procedures fulfill  requirements listed in Section~\ref{motivation}.
\subsection{Overview and notations}\label{notations}
 Our paper is organized as follows. In Section~\ref{explicationGLM}, we present the Goldenshluger and Lepski methodology in the setting of conditional density estimation. In Sections~\ref{kernel:section-} and \ref{sec:projection} respectively, kernel and projection rules are derived and studied in the oracle setting by using assumptions of Section \ref{assumptions:section}. Rates of convergence on anisotropic H\"{o}lder balls are studied in Section~\ref{minimax}. Then a simulation study is lead in Section~\ref{simulations}, where we focus on tuning aspects of our procedures. Finally, in Section \ref{sec:proofs} and in Appendix, we prove our results. To avoid too tedious technical aspects, most of proofs are only given for $d_1=d_2=1$ but can easily be extended to the general case. In the sequel, we assume that the sample size is $2n$. The first $n$ observations $(X_1,Y_1),\ldots,(X_n,Y_n)$ are used to estimate $f$, whereas $X_{n+1},\ldots,X_{2n}$ are used to estimate $\fx$ when necessary. We recall that for
any $i$, $X_i\in\R^{d_1}$ and $Y_i\in\R^{d_2}$ and we set $d=d_1+d_2$.

In addition to notations $R_x(\cdot,\cdot)$ and $\|\cdot\|_{x,2}$ introduced in  \eqref{riskx2} and \eqref{normx2}, we use for any $1\leq q<\infty$   $\|\cdot\|_q$,  the classical $\L_q$-norm of any function $g$: 
$$\|g\|^q_q=\int |g(x)|^qdx.$$
Some assumptions on functions $f$ and $\fx$, specified in Section~\ref{assumptions:section},  will depend on the following neighborhood of $x$, denoted $V_n(x)$: Given $A$ a positive real number and $(k_n)_n$ any positive sequence larger than $1$ only depending on $n$ and such that $k_n$ goes to $+\infty$, we set:
$$V_n(x)=\prod_{i=1}^{d_1}\left[x_i-\frac{2A}{k_n},x_i+\frac{2A}{k_n}\right].$$
Note that the size of $V_n(x)$ goes to 0. Then, we set
$$\|f\|_\infty=\sup_{t\in V_n(x)}\sup_{y\in\R^{d_2}}f(t,y)\in [0,+\infty], \quad \|\fx\|_\infty=\sup_{t\in V_n(x)}\fx(t)\in[0,+\infty]$$
and
$$\delta=\inf_{t\in V_n(x)}\fx(t)\geq 0.$$
Our results will strongly depend on these quantities. Finally, for any $u\in\R$, we set  $\{u\}_+=\max(u,0)$. 
\section{
Methodology}\label{explicationGLM}
\subsection{The Goldenshluger-Lepski methodology}\label{intro:methodo}
This section is devoted to the description of the Goldenshluger-Lepski methodology (GLM for short) in the setting of conditional density estimation.

The GLM consists in selecting an estimate from a family of estimates, each of them depending on a parameter $m$. Most of the time, choosing this tuning parameter can be associated to a regularization scheme: if we take $m$ too small, then the estimate oversmooths;  if we take $m$ too large, data are overfitted.

So, given a set of parameters $\mathcal{M}_n$, for any $m\in\mathcal{M}_n$, we assume we are given a smoothing linear operator denoted $\mathcal{K}_m$ and an  estimate $\hat f_m$. For any $m\in\mathcal{M}_n$, $\hat f_m$ is related to $\mathcal{K}_m(f)$ via its expectation and we assume that $\E[\hat f_m]$ is close to (or equal to) $\mathcal{K}_m(f)$.  The main assumptions needed for applying the GLM are
\begin{equation}\label{commute1}
\mathcal{K}_m\circ\mathcal{K}_{m'}=\mathcal{K}_{m'}\circ \mathcal{K}_m
\end{equation}
and
\begin{equation}\label{commute2}
\mathcal{K}_m(\hat f_{m'})=\mathcal{K}_{m'}(\hat f_m)
\end{equation}
for any $m,m'\in\mathcal{M}_n$.  The GLM  is a convenient way to select an estimate among $(\hat f_m)_{m\in\mathcal{M}_n}$ which amounts to selecting $m\in\mathcal{M}_n$ and can be described as follows:
For $\|\cdot\|$ a given norm and $\sigma$ a function to be chosen later, we set for any $m$ in $\mathcal{M}_n$,
$$A(m):=\sup_{m'\in\mathcal{M}_n}\left\{\|\hat f_{m'}-\mathcal{K}_{m'}(\hat f_m)\|-\sigma(m')\right\}_+.$$
Then we estimate $f$ by using $\hat f:=\hat f_{\hat m}$, where $\hat m$ is selected as follows:
$$\hat m:=\argmin_{m\in\mathcal{M}_n}\left\{A(m)+\sigma(m)\right\}.$$
This choice can be seen as a bias-variance tradeoff, with $\sigma(m)$ an estimator of the standard deviation of $\hat f_{m}$ and $A(m)$ an estimator of the bias (see later).
Let us now fix $m\in\mathcal{M}_{n}.$ Using (\ref{commute2}), we have:
\begin{eqnarray*}
\|\hat f-f\|&=&\|\hat f_{\hat m}-f\|\\
&\leq&\|\hat f_{\hat m}-\mathcal{K}_{\hat m}(\hat f_m)\|+\|\mathcal{K}_m(\hat f_{\hat m})-\hat f_m\|+\|\hat f_m-f\|\\
&\leq &A(m)+\sigma(\hat m)+A(\hat m)+\sigma(m)+\|\hat f_m-f\|\\
&\leq&2A(m)+2\sigma(m)+\|\hat f_m-\mathcal{K}_m(f)\|+\|\mathcal{K}_m(f)-f\|.
\end{eqnarray*}
But
\begin{eqnarray*}
A(m)&=&\sup_{m'\in\mathcal{M}_n}\left\{\|\hat f_{m'}-\mathcal{K}_{m'}(\hat f_m)\|-\sigma(m')\right\}_+\\
&\leq&\xi(m)+B(m)
\end{eqnarray*}
with  for any $m\in\mathcal{M}_n$,
$$\xi(m):=\sup_{m'\in\mathcal{M}_n}\left\{\|(\hat f_{m'}-\mathcal{K}_{m'}(f))-(\mathcal{K}_{m'}(\hat f_m)-(\mathcal{K}_{m'}\circ \mathcal{K}_m)(f))\|-\sigma(m')\right\}_+$$
and
$$B(m):=\sup_{m'\in\mathcal{M}_n}\|\mathcal{K}_{m'}(f)-(\mathcal{K}_{m'}\circ \mathcal{K}_m)(f)\|.$$
We finally obtain:
\begin{equation}\label{or1}
\|\hat f-f\|\leq 2B(m)+2\sigma(m)+\|\hat f_m-\mathcal{K}_m(f)\|+\|f-\mathcal{K}_m(f)\|+2\xi(m).
\end{equation}
Now, let us assume that
\begin{equation}\label{norme-op}
\vvvert\mathcal{K}\vvvert:=\sup_{m\in\mathcal{M}_n}\vvvert\mathcal{K}_m\vvvert<\infty,
\end{equation}
where $\vvvert\mathcal{K}_m\vvvert$ is the operator norm of  $\mathcal{K}_m$ associated with $\|\cdot\|.$ In this case, $B(m)$ is upper bounded by $\|f-\mathcal{K}_m(f)\|$ up to the constant $\vvvert\mathcal{K}\vvvert$, which corresponds to the bias of $\hat f_m$ if
\begin{equation}\label{E-S}
\mathcal{K}_m(f)=\E[\hat f_m].
\end{equation}
Furthermore, using (\ref{commute1}) and  (\ref{commute2}), for any $m\in\mathcal{M}_n$,
 $$\xi(m)\leq\sup_{m'\in\mathcal{M}_n}\left\{(1+\vvvert\mathcal{K}\vvvert)\|\hat f_{m'}-\mathcal{K}_{m'}(f)\|-\sigma(m')\right\}_+.$$
Then we choose $\sigma$ such that, with high probability, for any $m\in\mathcal{M}_n$,
\begin{equation}\label{concentration}
\|\hat f_m-\mathcal{K}_m(f)\|\leq \sigma(m)/(\vvvert \mathcal{K}\vvvert+1).
\end{equation}
So, (\ref{or1}) gives that, with high probability,
\begin{equation}\label{or2}
\|\hat f-f\|\leq C\inf_{m\in\mathcal{M}_n}\left\{\|f-\mathcal{K}_m(f)\|+\sigma(m)\right\},
\end{equation}
where $C$ depends only on $\vvvert\mathcal{K}\vvvert$. Since under (\ref{E-S}), $\sigma(m)$ controls the fluctuations of $\hat f_m$ around its expectation, $\sigma^{2}(m)$ can be viewed as a variance term and the oracle inequality (\ref{or2}) justifies our procedure. Previous computations combined with the upper bound of $A(m)$ also justify why $A(m)$ is viewed as an estimator of the bias.

Now, we illustrate this methodology with two natural smoothing linear operators: convolution and projection. The natural estimates associated with these operators are kernel rules and projection rules respectively. Next paragraphs  describe the main aspects of both procedures and discuss assumptions (\ref{commute1}),  (\ref{commute2}), (\ref{norme-op}) and (\ref{E-S}) which are the key steps of the GLM.

\subsection{Convolution and kernel rules}\label{convolution}
 Kernel rules are the most classical procedures for conditional density estimation. To estimate $f$, the natural approach consists in considering the ratio of a kernel estimate of $f_{X,Y}$ with a kernel estimate of $\fx$. Actually, we use an alternative approach and  to present our main ideas, we assume for a while that $\fx$ is known and positive.

  We introduce a kernel $K$, namely a bounded integrable function $K$ such that $\iint K(u,v)dudv=1$ and $\|K\|_2<\infty$.  Then, given a regularization parameter, namely a $d$-dimensional bandwidth $h$ belonging to a set ${\mathcal H}_n$ to be specified later, we set
  $$K_h(u,v)=\frac{1}{\prod_{i=1}^{d}h_i}K\left(\frac{u_1}{h_1},\ldots,\frac{u_{d_1}}{h_{d_1}}, \frac{v_1}{h_{d_1+1}},\ldots,\frac{v_{d_2}}{h_{d}}\right),\quad u\in\R^{d_1}, v \in\R^{d_2}.$$
 Then, we use the setting of Section \ref{intro:methodo} except that regularization parameters are denoted $h$, instead of $m$ to match with usual notation of the literature. Similarly, the set of bandwidths is denoted by ${\mathcal H}_n$, instead of ${\mathcal M}_n$. For any $h\in \mathcal{H}_n$, we set:
 $$\forall\, g\in\L_2,\quad \mathcal{K}_h(g)=K_h * g$$
 where  $*$ denotes the standard convolution product  and
\begin{equation}\label{def:methodo}
\hat f_h(x,y):=\frac{1}{n}\sum_{i=1}^n\frac{1}{\fx(X_i)}K_h(x-X_i,y-Y_i).
\end{equation}
The regularization operator $ \mathcal{K}_h$ corresponds to the convolution with $K_h$. Note that
$$\E[\hat{f}_{h}(x,y)]=(K_h * f)(x,y).$$
Therefore 3 of 4 assumptions of the GLM are satisfied, namely (\ref{commute1}), (\ref{commute2}) and (\ref{E-S}). Unfortunately, (\ref{norme-op}) is satisfied with $\|\cdot\|$ the classic $\L_2$-norm but not with $\|\cdot\|_{x,2}$, as adopted in this paper. We shall see how to overcome this problem later on.

Another drawback of this description is that $\hat f_h$ is based on the knowledge of $\fx$. A kernel rule based on $\hfx$, an estimate of $\fx$, is proposed in Section \ref{kernel:section} where we define $\sigma$ (see (\ref{concentration})) to apply the GLM methodology and then to obtain oracle inequalities similar to (\ref{or2}). Additional terms in oracle inequalities will be the price to pay for using $\hfx$ instead of $\fx$.

\subsection{Projection}\label{projection}
We introduce a collection of models $(S_m)_{m\in\mathcal{M}_n}$ and for any $m$, we denote $\mathcal{K}_m$ the projection on $(S_m,<,>_X)$ where $<,>_X$ is the scalar product defined by:
\begin{equation}\label{produitscalaire}
\forall g,g',\quad <g,g'>_X=\iint g(u,y)g'(u,y)\fx(u)dudy.
\end{equation}
Of course, (\ref{commute1}) is satisfied, but as for kernel rules, (\ref{norme-op}) is not valid with $\|\cdot\|=\|\cdot\|_{x,2}$.  Now, we introduce the following empirical contrast:
\begin{equation*}
 \text{ for all function } t,\quad \gamma_n(t)=
\frac{1}{n}\sum_{i=1}^{n}\left[\int_\mathbb{R} t^2(X_i,y)dy-2t(X_i, Y_i)\right],
\end{equation*}
so that $\E(\gamma_n(t))$ is minimum when $t=f$ (see Lemma~\ref{obvious} in Section  \ref{projection:section}).  Given $m$ in $\mathcal{M}_n$, the conditional density can be estimated by:
\begin{equation}\label{probmin}
\hat f_m \in\argmin_{t\in S_m}\gamma_n(t).
\end{equation}
Unlike kernel rules, this estimate does not depend on $\fx$ but (\ref{commute2}) and (\ref{E-S}) are not satisfied even if for large values of $n$, $\mathcal{K}_m(f)\approx\E[\hat f_m]$. Therefore, we modify this approach to overcome this problem.
The idea is the following. Let us denote  $S_{m\wedge m'}=S_{m}\cap S_{m'}$.
Taking inspiration from the fact that $\mathcal{K}_m\circ\mathcal{K}_{m'}(f)=\mathcal{K}_{m\wedge m'}(f)$, set for any $(m,m')\in{\mathcal M}_n^{2}$,
$$\tilde{\mathcal{K}}_m(\hat{f}_{m'})=\hat{f}_{m\wedge m'}.$$ This operator is only defined on the set of the estimators $\hat f_{m}$ but verifies
\eqref{commute2}.
Now the previous reasoning can be reproduced and the GLM described in Section \ref{intro:methodo} can be applied by replacing $\mathcal{K}_{m'}$ by $\tilde{\mathcal{K}}_{m'}$ in $A(m)$ and by setting
$$\xi(m):=\sup_{m'\in\mathcal{M}_n}\left\{\|(\hat f_{m'}-\mathcal{K}_{m'}(f))-(\tilde{\mathcal{K}}_{m'}(\hat f_m)-(\mathcal{K}_{m'}\circ \mathcal{K}_m)(f))\|-\sigma(m')\right\}_+.$$
In Section \ref{projection:section}, we define $\sigma$ such that for all $m,m'\in \mathcal{M}_{n}$, $\sigma(m \wedge m')\leq \sigma(m')$ and similarly to (\ref{concentration}), with high probability, for any $m\in\mathcal{M}_n$,
$$
\|\hat f_m-\mathcal{K}_m(f)\|\leq \frac{\sigma(m)}{2}.
$$
Then, for all $m,m'\in \mathcal{M}_{n}$,
$$\|\tilde{\mathcal{K}}_{m'}(\hat f_m)-(\mathcal{K}_{m'}\circ \mathcal{K}_m)(f)\|=\|\hat f_{m\wedge m'}-\mathcal{K}_{m\wedge m'}(f)\|\leq \frac{\sigma(m \wedge m')}{2}\leq \frac{\sigma(m')}{2}$$
so that $\xi(m)$ vanishes with high probability. Thus, we shall be able to derive oracle inequalities in this case as well.

\subsection{Discussion}
We have described two estimation schemes for which the GLM is appropriate: kernel and projection rules. In these schemes, the main commutative properties of the GLM, namely (\ref{commute1}) and (\ref{commute2}), are satisfied. Due to the particular choice of the loss-function $\|\cdot\|_{x,2}$, the property (\ref{norme-op})  is not satisfied. However in both schemes, we shall be able to prove that for any function $g$
\begin{equation}\label{K2>sup}
\| \mathcal{K}_m(g)\|_{x,2}\leq C\sup_{t\in V_n(x)} \|g\|_{t,2}
\end{equation}
where $C$ is a constant, $V_n(x)$ is the neighborhood of $x$ introduced in Section~\ref{notations}, and this property will allow us to control the bias term $B(m)$, as well as the term $\xi(m)$. In the sequel, we shall cope with the following specific features of each scheme:
\begin{itemize}
\item For kernel rules,  when $\fx$ is known, (\ref{E-S}) is satisfied and these estimates lead to straightforward application of the GLM. But, when $\fx$ is unknown, serious difficulties will arise.
\item For projection rules, the dependence on the knowledge of $\fx$ will be weaker but since (\ref{E-S}) is not satisfied, the control of the bias term will  not be straightforward.
\end{itemize}
Beyond these aspects, our main task in next sections will be to derive for each estimation scheme a function $\sigma$ that conveniently controls the fluctuations of preliminary estimates as explained in Section \ref{intro:methodo}.

\section{Assumptions}\label{assumptions:section}
 In this section, we state our assumptions
 on $f$ and $\fx$.
\begin{itemize}
\item[$(H_1)$] The conditional density $f$ is uniformly bounded on $V_n(x)\times\R^{d_2}$: $\|f\|_\infty<\infty.$
\item[$(H_2)$] The density $\fx$ is uniformly bounded on $V_n(x)$: $\|\fx\|_\infty<\infty.$
\item[$(H_3)$]  The density $\fx$ is bounded away from 0 on $V_n(x)$: $\delta>0.$ In the sequel, without loss of generality, we assume that $\delta\leq 1$.
\end{itemize}
Assumptions $(H_1)$ and $(H_2)$ are very mild. Note that under $(H_1)$, since $f$ is a conditional density, for any $t\in V_n(x)$, $\int f(t,v)dv=1$ and
\begin{equation}\label{sup>2}
\sup_{t\in V_n(x)}\|f\|_{t,2}^2\leq \sup_{t\in V_n(x), y\in\R^{d_2}}f(t,y)\int f(t,v)dv= \|f\|_\infty<\infty.
\end{equation}
 Assumption $(H_3)$ is not mild but is in some sense unavoidable. As said in Introduction, one goal of this paper is to measure the influence of the parameter $\delta$ on the performance of the estimators of $f$.

 For the procedures considered in this paper, if $\fx$ is unknown,  we need a preliminary estimator of $\fx$ denoted $\hfx$ that is constructed with observations $(X_i)_{i=n+1,\ldots, 2n}$. Then, we first assume that $\hfx$ satisfies the following condition:
\begin{equation}\label{defdeltachapeau}
\hat{\delta}:=\inf_{t\in V_n(x)}|\hfx(t)|>0.
\end{equation}
For estimating $\fx$, $\hfx$ has to be rather accurate:
\begin{equation}\label{Cfx}
\forall\,\lambda>0,\quad\P\left(\sup_{t\in V_n(x)}\left|\frac{\fx(t)-\hfx(t)}{\hfx(t)}\right|>\lambda\right)\leq \kappa\exp\{-(\log n)^{3/2}\},
\end{equation}
where $\kappa$ is a constant only depending on $\lambda$ and $\fx$.
Theorem~4 in \cite{nous} proves the existence of an estimate $\hfx$ satisfying these properties.

\section{Kernel rules}\label{kernel:section-}
In this section, we study the data-driven kernel rules we propose for estimating the conditional density $f$. They are precisely defined in Section \ref{kernel:section} and their theoretical performances in the oracle setting are studied in Section \ref{oraclekernel:section}. Before doing this, in Section \ref{lowerkernel:section}, we establish a lower bound of the risk for any kernel estimate.

\subsection{Lower bound for kernel rules}\label{lowerkernel:section}
In this section, we consider the kernel estimate $\hat f_h$ defined in (\ref{def:methodo}) for $h\in{\mathcal H}_n$. In particular, $\fx$ is assumed to be known. %
For any fixed $h\in\mathcal{H}_n$, we provide a lower bound of the risk of $\hat f_h$ with $q=2$ by using the following bias-variance decomposition:
$$R_x^2(\hat f_h,2)=\E\left[\|\hat f_h-f\|_{x,2}^2 \right]=\|K_h*f-f\|_{x,2}^2+\int\var(\hat{f}_{h}(x,y))dy.$$
\begin{proposition}\label{borneinforacle}
Assume that $(H_1)$ is satisfied. Then if $K(x,y)=K^{(1)}(x)K^{(2)}(y)$ with $K^{(1)}$ supported by $[-A,A]^{d_1}$, for any $h\in\mathcal{H}_n$, we have for any $n$,
$$R_x^2(\hat f_{h},2)\geq \|K_h*f-f\|_{x,2}^2+\frac{\|K^{(2)}\|_2^2}{n\prod_{i=1}^dh_i}\times \int\frac{[K^{(1)}(s)]^2}{\fx(x-(s_1h_1,\ldots,s_{d_1}h_{d_1}))}ds+\frac{C_1}{n},$$
where $C_1$ depends on $\|K\|_1$ and $\|f\|_{\infty}$.
If we further assume that $\fx$ is positive and continuous on a neighborhood of $x$, then if $\max\mathcal{H}_n\to 0$ when $n\to +\infty$,
\begin{equation}\label{minoopt}
R_x^2(\hat f_{h},2)\geq \|K_h*f-f\|_{x,2}^2+\frac{\|K\|_2^2}{\fx(x) n\prod_{i=1}^dh_i}\times (1+o(1))+O\left(\frac{1}{n}\right),
\end{equation}
when $n\to +\infty$.
\end{proposition}
The proof of Proposition~\ref{borneinforacle} is given in Section~\ref{sec:proofs}.
The lower bounds of Proposition~\ref{borneinforacle} can be viewed as benchmarks for our procedures. In particular, our challenge is to build a data-driven kernel procedure whose risk achieves the lower bound given in (\ref{minoopt}). It is the goal of the next section where we modify $\hat f_h$ by estimating $\fx$ when $\fx$ is unknown.

\subsection{Kernel estimator}\label{kernel:section}
Let us now define more precisely our kernel estimator. We consider the kernel $K$ defined in Section~\ref{convolution}, but following assumptions of Proposition~\ref{borneinforacle}, we further assume until the end of the paper that following conditions are satisfied.
\begin{itemize}
\item The kernel $K$ is of the form $K(u,v)=K^{(1)}(u)K^{(2)}(v),\quad u\in\R^{d_1}, v \in\R^{d_2}.$
\item The function $K^{(1)}$ is supported by $[-A,A]^{d_1}$.
\end{itemize}
Our data-driven procedure is based on $\hfx$ (see Section \ref{assumptions:section}) and is defined in the following way. We naturally replace $\hat{f}_h$ defined in (\ref{def:methodo}) with
\begin{equation}\label{estim1}
\hat{f}_h(x,y)=\frac{1}{n}\sum_{i=1}^{n}\frac{1}{\hfx(X_i)}K_{h}(x-X_i,y-Y_i).
\end{equation}
Then, we set
\begin{equation}\label{defsigma}\sigma(h)=\frac{\chi}{\sqrt{\hat{\delta}n\prod_{i=1}^dh_i}} \quad\mbox{with}\quad\chi=(1+\eta)(1+\|K\|_1)\|K\|_2,
\end{equation}
where $\hat\delta$ is defined in (\ref{defdeltachapeau}) and $\eta>0$ is a tuning parameter. The choice of this parameter will be discussed in Section~\ref{simulations} but all theoretical results are true for any $\eta>0$. We also specify the set $\mathcal{H}_n$:
\begin{itemize}
\item[$(CK)$] For any $h=(h_1,\ldots,h_d)\in \mathcal{H}_n$, we have for any $i$, $h_i^{-1}$ is a positive integer and
$$k_{n}\leq \frac{1}{h_i}, \ \forall\, � i\in \{1,\ldots,d_1\}, \quad\frac{1}{\prod_{i=1}^{d_1}h_i}\leq \frac{\hat{\delta} n}{(\log n)^{3}}\quad\mbox{and}\quad\log^2(n)\leq \frac{1}{\prod_{i=d_1+1}^{d}h_i}\leq n.$$
\end{itemize}
The GLM described in Section~\ref{convolution} can be applied and we estimate $f$ with $\hat{f}=\hat{f}_{\hat{h}}$ where
\begin{equation*}
\hat{h}=\hat{h}(x):=\argmin_{h\in\mathcal{H}_n}\left\{A(h)+\sigma(h)\right\},
\end{equation*}
\begin{equation*}
A(h):=\sup_{h'\in\mathcal{H}_n}\left\{\left\|\hat{f}_{h'}-\hat{f}_{h,h'}\right\|_{x,2}-\sigma(h')\right\}_+,
\end{equation*}
and
\begin{equation}\label{estim2}\hat{f}_{h,h'}(x,y)=\frac{1}{n}\sum_{i=1}^{n}\left[\hfx(X_i)\right]^{-1}(K_{h}*K_{h'})(x-X_i,y-Y_i)=(K_{h'}*\hat{f}_h)(x,y).
\end{equation}
In the case where $\fx$ is known, $\hfx$ is replaced by $\fx$ and $\hat\delta$ by $\delta$. In particular, we obtain the expressions of Section~\ref{convolution} except that now $\sigma$ is specified.

\subsection{Oracle inequalities for kernel rules}\label{oraclekernel:section}
We establish in this section oracle inequalities for our estimator $\hat f$ with in mind the benchmarks given in (\ref{minoopt}). 
To shed lights on the performance of our procedure and on the role of $\delta$, we first deal with the case where  $\fx$ is known. We first state a trajectorial oracle inequality and then a control of the risk.
\begin{theorem}\label{theo1}
Assume that the density $\fx$ is known so that $\hfx=\fx$. We also assume that $(H_1)$, $(H_3)$ and $(CK)$ are satisfied. If  $\delta n\geq 1$, we have with probability larger than $1-C\exp\{-(\log n)^{5/4}\},$
\begin{equation}
\|\hat{f}-f\|_{x,2}\le \inf_{h\in\mathcal{H}_n}\left\{C_1\sup_{t\in V_n(x)}\|K_h*f-f\|_{t,2}+\frac{C_2}{\sqrt{\delta n\prod_{i=1}^dh_i}}\right\},\label{traj1}
\end{equation}
where $C_1=1+2\|K\|_1$, $C_2=(1+\eta)\|K\|_2(3+2\|K\|_1)$ and $C$ depends on $K$, $\eta$ and $\|f\|_\infty$. Furthermore, for any $\pp\ge 1$,
\begin{equation}\label{majoopt}
R_x(\hat f,q) \le \tilde C_1\inf_{h\in\mathcal{H}_n}\left\{\sup_{t\in V_n(x)}\|K_h*f-f\|_{t,2}+\frac{1}{\sqrt{\delta n\prod_{i=1}^dh_i}}\right\}+ \frac{ \tilde C_2}{\sqrt{n}},
\end{equation}
where $ \tilde C_1$ depend on $K$, $\eta$ and $\pp$ and $ \tilde C_2$ depends on $K$, $\eta$, $\|f\|_\infty$ and $\pp$.
\end{theorem}
Due to the assumptions on ${\mathcal H}_n$, the last term of the right hand side of (\ref{majoopt}), namely $\tilde C_2/\sqrt{n}$, is negligible with respect to the first one.  Furthermore, since $\sigma^2(h)$ is proportional to $(\delta n\prod_{i=1}^nh_i)^{-1}$, the latter can be viewed as a variance term (see Section~\ref{intro:methodo}). Then right hand sides of (\ref{traj1}) and (\ref{majoopt}) correspond to the best tradeoff between a bias term and a variance term, so (\ref{traj1}) and (\ref{majoopt}) correspond indeed to oracle inequalities. Next, we can compare  the (squared) upper bound of (\ref{majoopt}) and the lower bound of (\ref{minoopt}) when $q=2$ and $\fx$ is continuous. We note that these bounds match up to leading constants, asymptotically negligible terms and up to the fact that terms of (\ref{majoopt}) are computed on $V_n(x)$ instead at $x$ (note that the size of $V_n(x)$ goes to 0 when $n\to +\infty$ and $\delta$ and $\fx(x)$ are close). Actually, since  (\ref{norme-op}) is not valid
for $\|\cdot\|=\|\cdot\|_{x,2}$, we
use Inequality (\ref{K2>sup}). This explains why we need to compute suprema of the bias term on $V_n(x)$.
Theorem \ref{theo1} shows the optimality of our kernel rule.

From these results, we can also draw interesting conclusions with respect to the term $\delta$ that appears in the variance term.  From (\ref{minoopt}), we already know that the term $\delta$ is unavoidable. Of course, the lower $\delta$ the worse the performance of $\hat f$. Actually, in the oracle context, our setting is (roughly speaking) equivalent to the classical setting where $\fx$ is lower bounded by an absolute constant (see \cite{brunelcomtelacour07} for instance), but with $\delta n$ observations to estimate $f$ instead of $n$. A similar remark will hold in the minimax framework of Section \ref{minimax}.

The following theorem deals with the general case where $\fx$ is unknown and estimated by~$\hfx$.
\begin{theorem}\label{theo2}
We assume that $(H_1)$, $(H_2)$, $(H_3)$, $(CK)$ (\ref{defdeltachapeau}) and (\ref{Cfx}) are satisfied. If $\delta n\geq 1$, we have with probability larger than $1-C\exp\{-(\log n)^{5/4}\},$
\begin{equation}
\|\hat{f}-f\|_{x,2}\le \inf_{h\in\mathcal{H}_n}\left\{C_1\sup_{t\in V_n(x)}\|K_h*f-f\|_{t,2}+\frac{C_2}{\sqrt{\hat{\delta}n\prod_{i=1}^dh_i}}\right\}+\frac{C_3}{\delta}\sup_{t\in V_n(x)}|\hfx(t)-\fx(t)|,\label{traj2}
\end{equation}
where $C_1=1+2\|K\|_1$, $C_2=(1+\eta)\|K\|_2(3+2\|K\|_1)$, $C_3$ depends on $K$, $\eta$ and $\|f\|_\infty$ and $C$ depends on $K$, $\eta$, $\fx$ and $\|f\|_\infty$.
Furthermore, for any $\pp\ge 1$,
\begin{equation*}
R_x(\hat f,q) \le \tilde C_1\inf_{h\in\mathcal{H}_n}\left\{\sup_{t\in V_n(x)}\|K_h*f-f\|_{t,2}+\frac{1}{\sqrt{\delta n\prod_{i=1}^dh_i}}\right\}+ \frac{ \tilde C_2}{\delta}\E^{\frac{1}{\pp}}\left(\sup_{t\in V_n(x)}|\hfx(t)-\fx(t)|^\pp\right)+\frac{ \tilde C_3}{\sqrt{n}},
\end{equation*}
where $ \tilde C_1$ depend on $K$, $\eta$ and $\pp$, $ \tilde C_2$ depends on $K$, $\eta$, $\pp$ and $\|f\|_\infty$ and $ \tilde C_3$ depends on $K$, $\eta$, $\fx$, $\|f\|_\infty$ and $\pp$.
\end{theorem}
The main difference between Theorems~\ref{theo2}  and \ref{theo1} lie in the terms involving $\sup_{t\in V_n(x)}|\hfx(t)-\fx(t)|$ in right hand sides. Of course, if $\fx$ is regular enough, we can build $\hfx$ so that this term is negligible. But in full generality, this unavoidable term due to the strong dependence of $\hat f_h$ on $\hfx$, may be cumbersome. Therefore, even if Theorem \ref{theo1} established the optimality of kernel rules in the case where $\fx$ is known, it seems reasonable to investigate other rules to  overcome this problem.

\section{Projection rules}\label{sec:projection}
Unlike previous kernel rules that strongly depend on the estimation of $\fx$, this section presents estimates  based on the least squares principle. The dependence on $\hfx$ is only expressed via the use of $\hat\delta$ and $\|\hfx\|_\infty:=\sup_{t\in V_n(x)}|\hfx(t)|$. For ease of presentation, we assume that $d_1=d_2=1$ but following results can be easily extended to the general case (see Section~\ref{proj:rates}).

\subsection{Models}
As previously, we are interested in the estimation of $f$ when the first variable is in the neighborhood of $x$, so we still use $V_n(x)$ defined in Section \ref{notations}. We introduce a collection of models $(S_m)_{m\in{\mathcal M}_n}.$
\begin{definition}\label{defSm}
Let ${\mathcal M}_n$ be a finite subset of $\{0,1,2,\ldots\}^2$. For
each $m=(m_1,m_2)\in {\mathcal M}_n$ and given two $\L_2(\R)$-orthonormal systems of bounded functions $(\varphi_j^m)_{j\in J_m}$ and $(\psi_k^m)_{k\in K_m}$, we set
$$F_{m_1}={\rm Span}(\varphi_j^m, \ j\in J_m),\quad H_{m_2}={\rm Span}(\psi_k^m, \ k\in K_m)$$
and the model $S_m$ is
$$S_m=F_{m_1}\otimes H_{m_2}=
\left\{t, \quad t(x,y)=\sum_{j\in J_m}\sum_{k\in K_m} {a}_{j,k}^m \varphi_j^m(x)\psi_k^m(y), \; a_{j,k}^m\in {\mathbb R}\right\}.$$
Finally, we denote
$$D_{m_1}=|J_m|\quad\mbox{and}\quad D_{m_2}=|K_m|,$$
respectively the dimension of $F_{m_1}$ and $H_{m_2}$.
\end{definition}
In this paper, we only focus on systems $(\varphi_j^m)_{j\in J_m}$ based on Legendre polynomials.
More precisely, the estimation interval $[x-2A,x+2A]$ is split into $2^{m_1}$ intervals of length $4A2^{-m_1}$: $$I_{l}=I_{l}^m=\left[x-2A+4A(l-1)2^{-m_1},x-2A+4Al2^{-m_1}\right)\qquad l=1,\dots, 2^{m_1}.$$
Then $J_m=\{(l,d), l=1,\dots, 2^{m_1}, d=0,\dots,r\}$, $D_{m_1}=(r+1)2^{m_1}$ and for any $u$,
$$\varphi_{j}^m(u)=\varphi_{l,d}^m(u)=\sqrt{\frac{2^{m_1}}{2A}}\sqrt{\frac{2d+1}{2}}P_{d}(T_l(u))\1_{I_{l}}(u)$$
where $P_{d}$ is the Legendre polynomial with degree $d$ on $[-1,1]$, 
and $T_l$ is the affine map which transforms $I_l$ into $[-1,1]$.

In the $y$-direction, we shall also take piecewise polynomials. In the sequel, we only use the following two assumptions : for all $m,m'\in \mathcal{M}_n$,
  $D_{m_2}\leq D_{m_2'}\Rightarrow H_{m_2}\subset H_{m_2'}$, and
there exists a positive real number $\phi_2$ such that for all $m\in\mathcal{M}_n$ for all $u\in\R$,
 $$\sum_{k\in K_m}(\psi_k^m)^2(u)\leq \phi_2D_{m_2}.$$
 Note that this assumption is also true for $F_{m_1}$. Indeed the spaces spanned by the $\varphi_j^m$'s are nested and, for all $u\in [x-2A,x+2A]$,
$$\sum_{l=1}^{2^{m_1}}\sum_{d=0}^r\varphi_{l,d}^m(u)^2\leq
\frac{2^{m_1}}{2A}\sum_{d=0}^r \frac{2d+1}2=\frac{2^{m_1}}{4A}(r+1)^2= \frac{r+1}{4A} D_{m_1} $$
using properties of the Legendre polynomials. Therefore, with $\phi_1=(r+1)/(4A)$, for any $u\in [x-2A,x+2A]$,
$$\sum_j(\varphi_j^m)^2(u)\leq \phi_1D_{m_1}.$$

\subsection{Projection estimator}\label{projection:section}
As in \citep{brunelcomtelacour07} and following Section \ref{projection},  we introduce the following empirical contrast:
\begin{equation*}\label{gam0} \gamma_n(t)=
\frac{1}{n}\sum_{i=1}^{n}\left[\int_\mathbb{R} t^2(X_i,y)dy-2t(X_i, Y_i)\right].
\end{equation*}
We have the following lemma whose proof is easy by using straightforward computations. We use the norm $\|\cdot\|_{X}$ associated with the dot product $ \langle , \rangle_X$ defined in (\ref{produitscalaire}), so we have for any $t$, $$\|t\|_{X}^{2}=\iint t^2(u,y)\fx(u)dudy.$$
\begin{lemma}\label{obvious}
Assume that the function
$\sum_{j\in J_m}\sum_{k\in K_m}\hat {a}_{j,k}^m \varphi_j^m\psi_k^m$ minimizes the empirical contrast function $\gamma_n$ on $S_m$, then
\begin{equation}
\hat{G}_m\hat{A}_m= \hat{Z}_m,\label{derivation}
\end{equation}
where
$ \hat{A}_m$ denotes the matrix with coefficients $(\hat{a}_{j,k}^m)_{j\in J_m, k\in K_m}$,
$$\hat{G}_m=\left(\displaystyle
\frac{1}{n}\sum_{i=1}^n\varphi_{j_1}^m(X_i)\varphi_{j_2}^m(X_i)\right)
_{j_1, j_2 \in J_m}\quad \mbox{and}\quad
\hat{Z}_m=\left(\displaystyle\cfrac{1}{n}\sum_{i=1}^n\varphi_j^m(X_i)
\psi_k^m(Y_i)\right) _{j\in J_m, k \in K_m}. $$
Similarly, if $\mathcal{K}_m(f)$ is the orthogonal projection of $f$ on $(S_m, \langle , \rangle_X)$, it minimizes on $S_m$
$$t\longmapsto\gamma(t)=\|t-f\|_{X}^{2}-\|f\|_{X}^{2}=\E(\gamma_{n}(t))$$
and if  $\mathcal{K}_m(f)=\sum_{j\in J_m}\sum_{k\in K_m}a_{j,k}^m \varphi_j^m\psi_k^m$ then, $$G_mA_m=Z_m,$$ where $A_m$ denotes the matrix with coefficients $(a_{j,k}^m)_{j\in J_m, k\in K_m}$,
 $G_m=\E(\hat{G}_m)=\left(\langle
\varphi_{j_1}^m,\varphi_{j_2}^m
\rangle_{X}\right)_{j_1, j_2 \in J_m}$ and 
$$Z_m=\E(\hat{Z}_m)=\left(\iint
\varphi_{j}^m(u)\psi_k^m(y)
f(u,y)\fx(u)dudy\right) _{j\in J_m, k \in K_m}.$$
\end{lemma}
From this lemma, we obtain that $\E(\gamma_n(t))$ is minimum when $t=f$, which justifies the use of $\gamma_n$.

Then, to derive $\hat f_m$ an estimate of $f$, we use (\ref{derivation}) as a natural consequence of the minimization problem (\ref{probmin}). But if $\hat G_m$ is not invertible, $\hat A_m$ can be not uniquely defined.

Since $x$ is fixed, we can define, for each $m=(m_1,m_2)$, the index  $l_{m_1}=l_{m_1}(x)$ such that $x$ belongs to $I_{l_{m_1}}$ (actually, since the estimation interval is centered in $x$, $l_{m_1}=2^{m_1-1}+1$).
Furthermore, since we use a piecewise polynomial system, the Gram matrix $\hat{G}_m$ is a block diagonal matrix with blocks
 $\hat{G}_m^{(1)}, \dots, \hat{G}_m^{(2^{m_1})}$, where
$$\hat{G}_m^{(l)}=\left(\displaystyle
\frac{1}{n}\sum_{i=1}^n\varphi_{l,d_1}^m(X_i)\varphi_{l,d_2}^m(X_i)\right)
_{0\leq d_1,d_2 \leq r}.$$
In the same way, we can define for $l=1,\dots, 2^{m_1}$
$$\hat{Z}_m^{(l)}=\left(\displaystyle\cfrac{1}{n}\sum_{i=1}^n\varphi_{l,d}^m(X_i)
\psi_k^m(Y_i)\right) _{0\leq d\leq r, k \in K_m}.$$
Now,  and by naturally using the blockwise representation of $\hat G_m$, we define the collection of estimators $(\hat{f}_{m})_{m\in{\mathcal M}_n}$ as:
$$\hat{f}_{m}(x,y)=\sum_{d=0}^{r}\sum_{k\in K_m}\hat a_{(l_{m_1},d),k}^m\varphi_{l_{m_1},d}^m(x)\psi_{k}^{m}(y)$$
and $$(\hat a_{(l_{m_1},d),k}^m)_{0\leq d\leq r,k\in K_m}:=\hat A_{m}^{(l_{m_1})}:=\begin{cases}
(\hat{G}_m^{(l_{m_1})})^{-1}\hat{Z}_m^{(l_{m_1})}& \text{if } \min({\rm Sp}(\hat{G}_m^{(l_{m_1})}))>(1+\eta)^{-2/5} \hat{\delta}\\
0 & \text{otherwise,}
\end{cases}$$
where $\eta$ is a positive real number. Here, for a  symmetric matrix $M$, ${\rm Sp}(M)$ denotes the spectrum of $M$, i.e. the set of its eigenvalues. This expression allows us to overcome problems if $\hat G_m$ is not invertible.
Note that, when $r=0$, where $r$ is maximal degree of Legendre polynomials, this estimator can be written
 $$\hat f_m(x,y)=\sum_{j\in J_m}\sum_{k\in K_m}\cfrac{\sum_{i=1}^n\varphi_j^m(X_i)
\psi_k^m(Y_i)}{\sum_{i=1}^n\varphi_j^m(X_i)^2}
 \varphi_j^m(x)\psi_k^m(y).$$
Now, to choose a final estimator among this collection, as explained in Section \ref{projection}, we denote
$m\wedge j=(m_1\wedge j_1,m_2\wedge j_2)=(\min(m_1,j_1),\min(m_2,j_2))$ and by using $\hfx$ introduced in Section \ref{assumptions:section}, we set
\begin{equation}\label{chi}
\sigma(m)=\hat{\chi}\sqrt{\frac{D_{m_1}D_{m_2}}{ \hat{\delta} n}}\quad
\mbox{with}\quad \hat{\chi}^{2}=
(1+\eta)^2(4\phi_1\phi_2(r+1))\frac{\widehat{\|\fx\|_\infty}}{\hat{\delta}},
\end{equation}
where $\widehat{\|\fx\|_\infty}=\|\hfx\|_\infty$ and $\hat\delta$ is defined in \eqref{defdeltachapeau}.
We also specify the models we use: The following condition is the analog of $(CK)$:
\begin{enumerate}
\item[$(CM)$]
For any $m\in \mathcal{M}_n$,
$$k_{n}(r+1)\leq D_{m_1}\leq \frac{\hat{\delta} n}{(\log n)^{3}}\quad\mbox{and}\quad\log^2(n)\leq D_{m_2}\leq n.$$
\end{enumerate}
The GLM described in Section~\ref{convolution} can be applied and we estimate $f$ with $\tilde f=\hat{f}_{\hat{m}}$ where
\begin{equation*}
\hat{m}=\hat{m}(x):=\arg\min_{m\in\mathcal{M}_n}\left\{ A(m)+\sigma(m)\right\}
\end{equation*}
and
\begin{equation*}
A(m):=\sup_{m'\in \mathcal{M}_n}\left[\|\hat f_{m'}-\hat f_{m'\wedge m}\|_{x,2}-\sigma(m')\right]_+.
\end{equation*}
The next section studies the performance of the estimate $\tilde f$.

\subsection{Oracle inequality for projection estimators}
We establish in this section oracle inequalities for the projection estimate in the same spirit as for the kernel rule.
We recall that  $\mathcal{K}_m(f)$ is  the orthogonal projection of $f$ on $(S_m,<,>_X)$ where $<,>_X$ is the dot product defined in \eqref{produitscalaire}. The following result is the analog of Theorem~\ref{theo2}.
 \begin{theorem} \label{oracleproba}
 We assume that $(H_1)$, $(H_2)$, $(H_3)$, $(CM)$ (\ref{defdeltachapeau}) and (\ref{Cfx}) are satisfied. If $\delta n\geq 1$, we have with probability larger than $1-C\exp\{-(\log n)^{5/4}\},$
 \begin{equation*}
 \|\tilde{f}-f\|_{x,2}\leq \inf_{m\in\mathcal{M}_n}\left(C_1\sup_{t\in V_n(x)}\|\mathcal{K}_m(f)-f\|_{t,2}+\frac{5}{2}\hat{\chi}\sqrt{\frac{D_{m_1}D_{m_2}}{ \hat{\delta} n}}\right)
 \end{equation*}
 with $\hat{\chi}$  defined in (\ref{chi}), $C_1=1+2(r+1)\delta^{-1}\|\fx\|_\infty$ and $C$ depends on $\phi_1,\phi_2, r, \eta, \|f\|_\infty$ and $\fx$. Furthermore, for any $q\geq 1$
\begin{equation*}\label{risk}
R_x(\tilde{f},q)\leq \tilde C_1\inf_{m\in\mathcal{M}_n}\left(\sup_{t\in V_n(x)}\|\mathcal{K}_m(f)-f\|_{t,2}+\sqrt{\frac{D_{m_1}D_{m_2}}{\delta n}}\right)+\frac{\tilde C_2}{\sqrt{ n}}
\end{equation*}
where $\tilde C_1$ depends on $\phi_1,\phi_2,r, \eta, \|\fx\|_\infty, \delta$ and $q$ and $\tilde C_2$ depends on $\phi_1,\phi_2, r, \eta, \|f\|_\infty, \fx$ and $q$.
\end{theorem}
As for Theorem~\ref{theo2}, using the definition of $\sigma$, the right hand sides correspond to the best tradeoff between a bias term and a variance term. Note that unlike kernel rules,  the performances of $\tilde f$ do not depend on the rate of convergence of $\hfx$ for estimating $\fx$. But there is a price to pay:  due to a rougher control of the bias term, $\hat\chi$ depends on $\hat\delta$ and the leading constants $C_1$ and $\tilde C_1$ depend on $\delta$. In particular, when $\fx$ is known, conclusions drawn from Theorem~\ref{theo1} do not hold here. However, in the case where $r=0$ (the basis in the first coordinate is simply the histogram basis), we can use the simpler penalty term
$\hat{\chi}=(1+\eta)\sqrt{4\phi_1\phi_2}$ and the previous result still holds.
To prove this, it is sufficient to use the basis $(\|\varphi_j\|_X^{-1}\varphi_j\otimes \psi_k)_{j,k}$ which is orthonormal for the scalar product $\langle .,.  \rangle_X$.

\section{Rates of convergence}\label{minimax}
In this section, minimax rates of convergence will be computed on H\"{o}lder balls $\mathcal{H}_d(\boldsymbol{\alpha},\bold L)$. We recall that for
 two $d$-tuples of positive reals $\boldsymbol{\alpha}=(\alpha_1,\dots,\alpha_d)$ and $\bold L=(L_1,\dots,L_d)$,
 \begin{eqnarray*}
&&\mathcal{H}_d(\boldsymbol{\alpha},\bold L)=\Big\{ f: \R^d\to \R \text{  s.t. } \forall\, 1\leq i \leq d\quad
\left\|\frac{\partial^{m}f}{\partial x_i^{m}}\right\|_\infty\leq L_i, \quad m=0,\dots,\lfloor \alpha_i\rfloor\\
&& \text{ and for all } \,t\in \R \quad \left\|\frac{\partial^{\lfloor \alpha_i\rfloor}f}{\partial x_i^{\lfloor \alpha_i\rfloor}}(\cdot+te_i) - \frac{\partial^{\lfloor \alpha_i\rfloor}f}{\partial x_i^{\lfloor \alpha_i\rfloor}}(\cdot)\right\|_{\infty}\leq L_i |t|^{\alpha_i- \lfloor \alpha_i\rfloor}
\Big\}\\
\end{eqnarray*}
where for any $i$, $\lfloor \alpha_i\rfloor=\max\{l\in\mathbb{N}: l<\alpha_i\}$ and $e_i$ is the vector where all coordinates are null except the $i$th one which is equal to 1.
In the sequel, we use the classical anisotropic smoothness index defined by
$$\bar\alpha=\left(\sum_{i=1}^d\frac1{\alpha_i}\right)^{-1}$$
and introduced in the seminal paper \cite{KLP}. See also \cite{GL4bis}. 

\subsection{Lower bound}\label{seclow}
We have the following result that holds without making any assumption. It is proved in Section~8.3 of \cite{nous}.
\begin{theorem}\label{t2} There exists a positive constant $C$ not depending on $\bold{L}$  nor $n$ such that, if $n$ is large enough,
$$\inf_{T_n}\sup_{(f,\fx)\in \mathcal{\tilde H}(\boldsymbol{\alpha},\bold L)}\left\{(\fx(x))^{\frac{2\bar\alpha}{2\bar\alpha+1}}\E\|f-T_n\|_{x,2}^2\right\}\geq C\left(\prod_{i=1}^dL_i^{\frac{1}{\alpha_i}}\right)^{\frac{2\bar\alpha}{2\bar\alpha+1}}n^{-\frac{2\bar\alpha}{2\bar\alpha+1}},$$
where the infimum is taken over all estimators $T_n$ of $f$ based on the observations
$(X_i,Y_i)_{i=1,\ldots,n}$ and
 $\mathcal{\tilde H}(\boldsymbol{\alpha},\bold L) $ is the set such that the conditional density $f$ belongs to
 $\mathcal{H}_d(\boldsymbol{\alpha},\bold L)$ and the marginal density $\fx$
 is continuous.
 \end{theorem}
 Note that we consider the ball $\mathcal{\tilde H}(\boldsymbol{\alpha},\bold L)$ which may be (slightly) smaller than the ball $\mathcal{H}(\boldsymbol{\alpha},\bold L)$. Actually, we wish to point out the dependence of the lower bound with respect to $n$, $\boldsymbol{\alpha}$ and $\bold L$ as usual but also to $\fx(x)$, which is less classical. The goal in next sections is to show that our procedures achieve the lower bound of Theorem~\ref{t2}.

\subsection{Upper Bound for kernel rules}
In this section, we need an additional assumption on $f$.
\begin{itemize}
\item[$(H_4)$] There exists a compact set $B$, such that for all $t\in V_n(x)$, the function $y\mapsto f(t,y)$ has a support included into $B$. We denote by $|B|$ the length of the compact set $B$.
\end{itemize}
This assumption could be avoided at the price of studying the risk restricted on $B$.
Moreover, to study the bias of the kernel estimator, we consider for any $\bold{M}=(M_1,\ldots,M_d)$ the following condition.
\begin{itemize}
\item[$(BK_{\bold M})$] For  any $i\in\{1,\ldots,d\}$, for any $1\leq j\leq M_i$,
we have
\begin{equation*}\label{eq2}
\int_{\R} |x_i|^j|K(x)|dx_i<\infty\quad\mbox{and}\quad\int_{\R} x_i^jK(x)dx_i=0.
\end{equation*}
\end{itemize}
We refer the reader to \cite{KLP} for the construction of a kernel $K$ satisfying $(BK_{\bold M})$ and previous required conditions.
We obtain the following result showing the optimality of our first procedure from the minimax point of view, up to the rate for estimating $\fx$.
\begin{theorem}\label{propvitesse}
We assume that $(H_1)$, $(H_2)$, $(H_3)$, $(H_4)$, $(CK)$, (\ref{defdeltachapeau}) and (\ref{Cfx}) are satisfied. Let $\bold{M} =(M_1,\ldots,M_d)$ such that $(BK_{\bold M})$ is satisfied. Then if $f$ belongs to $\mathcal H_{d}(\boldsymbol\alpha,\bold L)$  such that $\lfloor\alpha_i\rfloor\le M_i$ for all $i=1,\ldots,d$, the kernel rule $\hat{f}$  satisfies
 for any $\pp\ge 1$,
\begin{equation*}
R_x^q(\hat f,q) \le \tilde C_1\left(\prod_{i=1}^dL_i^{\frac{1}{\alpha_i}}\right)^{\frac{q\bar\alpha}{2\bar\alpha+1}}(n\delta)^{-\frac{\pp\overline{\alpha}}{2\overline{\alpha}+1}}+ \frac{ \tilde C_2}{\delta^q}\E\left(\sup_{t\in V_n(x)}|\hfx(t)-\fx(t)|^\pp\right)+\tilde C_3n^{-\frac{q}{2}},
\end{equation*}
where $ \tilde C_1$ depend on $K$, $\eta$ and $\pp$, $ \tilde C_2$ depends on $K$, $\eta$, $\pp$ and $\|f\|_\infty$ and $ \tilde C_3$ depends on $K$, $\eta$, $\fx$, $\|f\|_\infty$ and $\pp$.
\end{theorem}
If the leading term in the last expression is the first one, then, up to some constants, the upper bound of Theorem~\ref{propvitesse} matches with the lower bound obtained in Theorem~\ref{t2} (note that $\delta$ is close to $\fx(x)$) when $q=2$. In this case, our estimate is adaptive minimax. To study the second term, we can use Theorem~4 of \cite{nous} that proves that, in our setting, there exists an estimate $\hfx$ achieving the rate $(\log /n)^{\frac{\bar\beta}{2\bar\beta+1}}$ if $\fx\in\mathcal{H}_{d_1}(\boldsymbol{\beta},\bold{\tilde L})$ and we obtain the following corollary.
\begin{cor}\label{cor}
We assume that $(H_1)$, $(H_2)$, $(H_3)$, $(H_4)$, $(CK)$ and $(BK_M)$ are satisfied. We also assume that  $\fx\in\mathcal{H}_{d_1}(\boldsymbol{\beta},\bold{\tilde L})$ such that for any $i=1,\ldots,d_1$, $\tilde L_i>0$ and
$0< \beta_i\le \beta^{(m)}_i$  with some known $\beta^{(m)}_i>0$.
 Then if $f$ belongs to $\mathcal H_{d}(\boldsymbol\alpha,\bold L)$  such that $\lfloor\alpha_i\rfloor\le M_i$ for all $i=1,\ldots,d$, the kernel rule $\hat{f}$  satisfies
 for any $\pp\ge 1$,
$$R_x^q(\hat f,q) \le C_1\left(\left(\prod_{i=1}^dL_i^{\frac{1}{\alpha_i}}\right)^{\frac{q\bar\alpha}{2\bar\alpha+1}}(n\delta)^{-\frac{\pp\overline{\alpha}}{2\overline{\alpha}+1}}+ \frac{1}{\delta^\pp}\left(\frac{\log n}{n}\right)^{\frac{\pp\bar\beta}{2\bar\beta+1}}\right)+C_2n^{-\frac{q}{2}},$$
where
$C_1$ is a constant not depending on $\bold L$, $n$ and $\delta$ and $C_2$ is a constant not depending on $\bold L$ and $n$.
\end{cor}
From the corollary, we deduce that if $\bar\beta>\bar\alpha$ and if $\delta$ is viewed as a constant, then the leading term is the first one. Furthermore, in this case, the rate is polynomial and the rate exponent is the classical ratio associated with anisotropic H\"{o}lder balls: $\bar\alpha/(2\bar\alpha+1)$. Our result also explicits the dependence of the rate with respect to $\bold{L}$ and $\delta$.

\subsection{Upper bound for projection estimates}\label{proj:rates}
In the same way, we can control the bias for our second procedure of estimation in order to study the rate of convergence.
Let us  briefly explain how the procedure defined in Section~\ref{sec:projection} can be extended to the estimation of conditional anisotropic densities $f:\R^{d_1}\times\R^{d_2}\to \R$ with $d_1,d_2\geq 2$. The contrast is still the same and
the estimators $\hat{f}_m$ have to be defined for $m=(m_1,\ldots,m_{d})$ with
a polynomial basis on hyperrectangles : see \cite{akakpolacour} for a precise definition. The model dimension is now
$$D_{m_1}=\prod_{i=1}^{d_1}r_i2^{m_i}$$
where $r_1,\dots, r_{d_1}$ are the maximum degrees.
Then, the selection rule to define $\tilde f$ is unchanged, except that in \eqref{chi}
$$\hat{\chi}^{2}=(1+\eta)^2\left(4\phi_1\phi_2\prod_{i=1}^{d_1}r_i\right)\frac{\widehat{\|\fx\|_\infty}}{\hat{\delta}}$$
In order to control precisely the bias, we introduce the following condition.
\begin{itemize}
\item[$(BM_{\boldsymbol r})$] $H_{m_{2}}$ is a space of piecewise polynomials with degrees bounded by $r_{d_1+1},\dots, r_d$, with
$D_{m_2}=\prod_{i=d_1+1}^{d}r_i2^{m_i}$.
\end{itemize}
This allows us to state the following result.
\begin{theorem}\label{vitesse}
We assume that $(H_1)$, $(H_2)$, $(H_3)$, $(H_4)$, $(CM)$, (\ref{defdeltachapeau}) and (\ref{Cfx}) are satisfied. Let $\boldsymbol{r} =(r_1,\ldots,r_d)$ such that $(BM_{\boldsymbol r})$ is satisfied. Then if $f$ belongs to $\mathcal H_{d}(\boldsymbol\alpha,\bold L)$  such that $\alpha_i< r_i$ for all $i=1,\ldots,d$, the projection rule $\tilde{f}$  satisfies
 for any $\pp\ge 1$,
\begin{equation*}
R_x^q(\tilde f,q) \le \tilde C\left(\prod_{i=1}^dL_i^{\frac{1}{\alpha_i}}\right)^{\frac{q\bar\alpha}{2\bar\alpha+1}}n^{-\frac{\pp\overline{\alpha}}{2\overline{\alpha}+1}},
\end{equation*}
where $ \tilde C$ depend on $A,|B|, \boldsymbol r,\boldsymbol\alpha,\delta$ and $\|\fx\|_\infty$.
\end{theorem}
Thus, even if the control of $\delta$ is less accurate, the projection estimator achieves the optimal rate of convergence
whatever the regularity of $\fx$.

\section{Simulations}\label{simulations}

In this section we focus on the numerical performances of our estimators. We  first describe the algorithms. Then, we introduce the studied examples and we illustrate the performances of our procedures with some figures and tables.

\subsection{Estimation algorithms}
For both methods (kernel or projection), we need a preliminary estimator of $\fx$. In order to obtain an accurate estimator
of $\fx$, we use a pointwise Goldenshluger Lepski procedure which consists in the following for estimating $\fx$ at $x$. This preliminary estimator is constructed using the sample $(X_i)_{i=n+1,\ldots,2n}$.
Let us define for $h>0$,
\begin{equation}
  \mathrm{pen}(n,h)= 2.2\|K\|_2(1+\|K\|_1)\sqrt{\frac{|\log h|\tilde{\mathrm{f}}_X(x)}{nh}},
\end{equation}
where $\tilde{\mathrm{f}}_X$ is a preliminary estimator of $\mathrm{f}_X$ obtained by the rule of thumb (see \cite{silverman}), and $K$ is the classical Gaussian kernel. The value 2.2 is the adjusted tuning constant which was convenient on a set of preliminary simulations.
 Given $H$ a finite set of bandwidths (actually $H$ is a set of $10$ bandwidths centered at the bandwidth obtained by the rule of thumb) and for $h,h'\in H$, consider
\begin{equation*}
  \hat {\mathrm{f}}_{h}(x)=\frac{1}{n}\sum_{i=n+1}^{2n}K_h(x-X_i)
\quad\text{ and }\quad
 \hat {\mathrm{f}}_{h,h'}(x) = \frac{1}{n} \sum_{i=n+1}^{2n}(K_h*K_{h'})(x-X_i).
\end{equation*}
 We consider
\begin{equation*}
\mathrm{A}(h,x):=\max_{h'\in H}\left\{\left|\hat{\mathrm{f}}_{h,h'}(x)-\hat{\mathrm{f}}_{h'}(x)\right|-\mathrm{pen}(n,h')\right\}_+.
\end{equation*}
Finally we define $h_0$ by
\begin{equation}\label{eq:equilibre}
h_0:=\argmin_{h\in H}\left\{\mathrm{A}(h,x)+\mathrm{pen}(n,h)\right\}
\end{equation}
and we consider the following procedure of estimation:
 $ \hfx(x) = \hat {\mathrm{f}}_{ h_0}(x).$

Now, the algorithm for the kernel estimation of $f$ is entirely described in Section~\ref{kernel:section} and we perform it with $K$ the Gaussian kernel and a set of
$10$ bandwidths in each direction, that means that the size of $\mathcal{H}_n$ is $10^{d_1+d_2}$. The quantity $\|\hat{f}_{h'}-\hat{f}_{h,h'}\|_{x,2}$ is made easy to compute with some preliminary theoretical computations (in particular, note that for the Gaussian kernel $K_h*K_h'=K_{h''}$ with $h''^2=h^2+h'^2$).
The only remaining parameter to tune is $\eta$ which appears in the penalty term $\sigma$ (see \eqref{defsigma}).

In the same way, we follow Section~\ref {projection:section} to implement the projection estimator. Matrix computations are easy to implement and make the implementation very fast. We only present the case of polynomials with degrees $r=s=0$, i.e. histograms, since the performance is already good in this case.  Again, the only remaining parameter to tune is $\eta$ which appears in the penalty term $\sigma$ (see \eqref{chi}). Note that in the programs, it is possible to use non-integers $m_i$ and in fact this improves the performance of the estimation. However, to match with the theory we shall not tackle this issue.

\subsection{Simulation study and analysis}

\begin{figure}[h]
\begin{tabular}{ll}
 \includegraphics[scale=0.37]{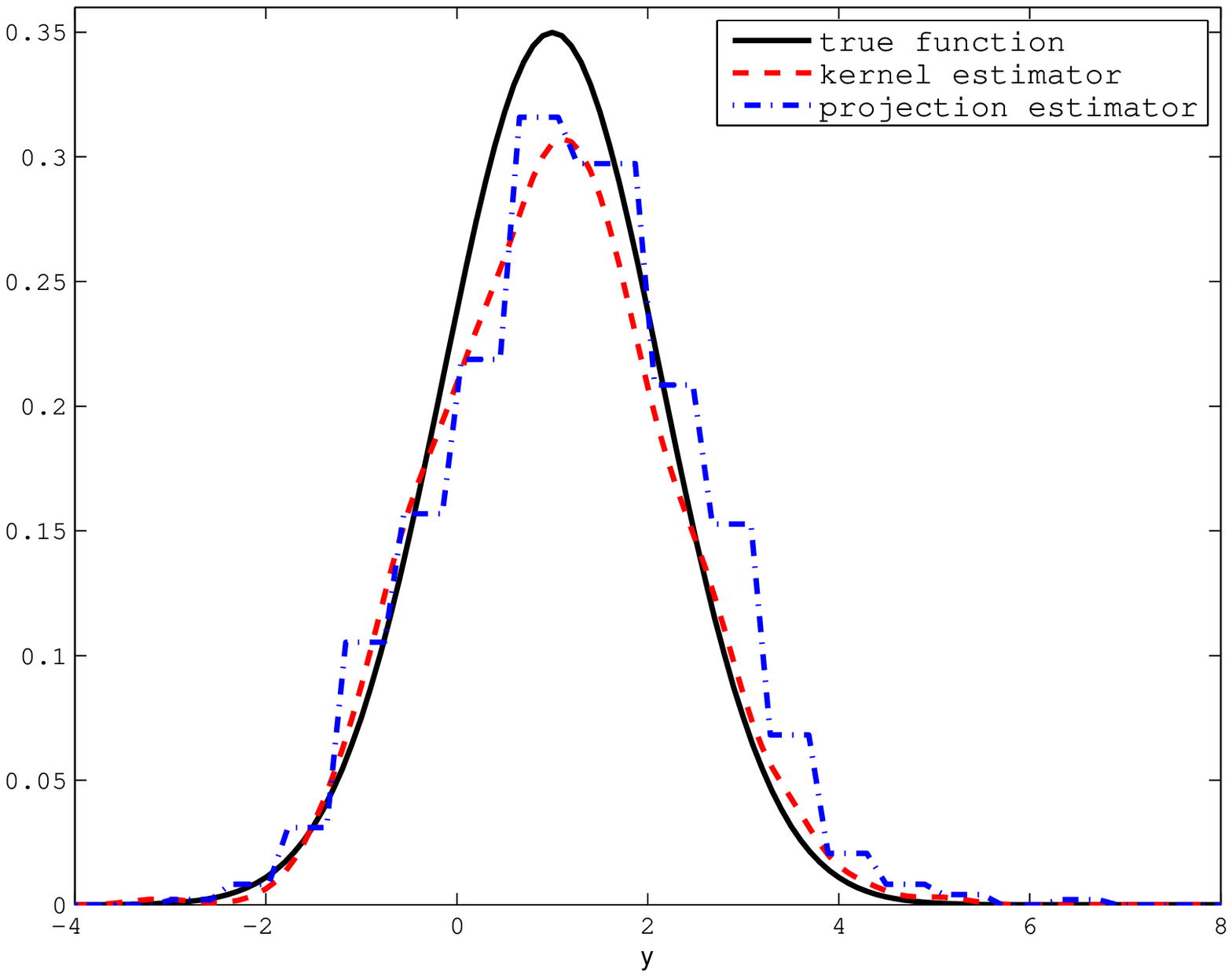}
& \includegraphics[scale=0.37]{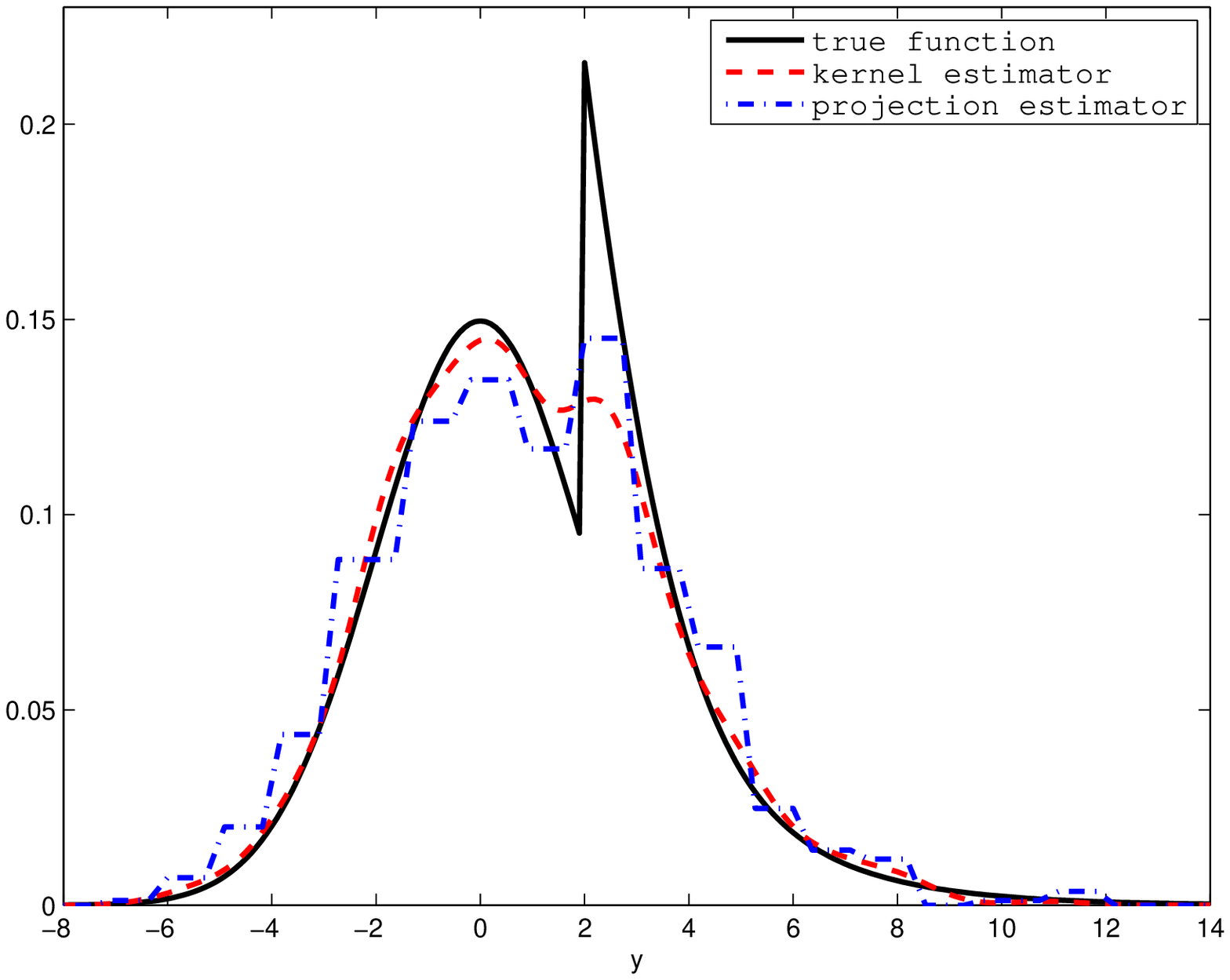}
\end{tabular}
 \caption{Plots of true function $f(x,.)$ (plain line) versus kernel estimator $\hat f(x,.)$ (dashed line)
 and projection estimator $\tilde f(x,.)$ (dot-dashed line) in $x=0$ ($n=1000$) for Example 3 (left) and Example 4 (right)}
 \label{reconstruction}
\end{figure}

We apply our procedures to different examples of conditional density functions with $d_1=d_2=1$.
More precisely, we observe $(X_i,Y_i)_{i=1,\ldots,n}$ such that
\begin{itemize}
\item[Example 1] The $X_i$'s are iid uniform variables on $[0,1]$ and
$$Y_i=2X_i^2+5+\e_i(1.3-|X_i|)^{1/2},\quad i=1,\ldots,n,$$
where the $\e_i$'s are i.i.d. reduced and centered Gaussian variables, independent of the
$X_i$'s.
Note that we also studied heavy-tailed noises in this example (i.e. the $\e_i$'s are variables with a standard Cauchy distribution) and the results were almost identical.

\item[Example 2] The $X_i$'s are iid uniform variables on $[0,1]$ and the distribution of the $Y_i$'s is a mixture of a normal distribution and an exponential distribution:
$Y_i\sim 0.75 \e_i+0.25 (2+E_i)$, where
$\e_i$ is a zero-mean normal distribution with standard deviation $2+X_i$ and $E_i$ is exponential with parameter $2$.
\item[Example 3] The $X_i$'s are iid and their common distribution is a mixture of two normal distributions,
$0.5\mathcal{N}(0,1/81)+0.5\mathcal{N}(1,1/16)$
 and
$$Y_i=X_i^2+1+\e_i(1.3+|X_i|)^{1/2},\quad i=1,\ldots,n,$$
where the $\e_i$'s are i.i.d. reduced and centered Gaussian variables, independent of the
$X_i$'s.
 \item[Example 4] The $X_i$'s are iid and their common distribution is a mixture of two normal distributions,
 $0.5\mathcal{N}(0,1/81)+0.5\mathcal{N}(1,1/16)$
 and the distribution of the $Y_i$'s is a mixture of a normal distribution and an exponential distribution:
$Y_i\sim 0.75 \e_i+0.25 (2+E_i)$, where
$\e_i$ is a zero-mean normal distribution with standard deviation $2+X_i$ and $E_i$ is exponential with parameter $2$.
\end{itemize}
We simulate our observations for three sample sizes: $n=250$, $n=500$ and $n=1000$.  In Figure~\ref{reconstruction}, we illustrate the quality of reconstructions for both estimates when $\fx$ is unknown.
We use $\eta=-0.2$ for the projection estimator and $\eta=1$ for the kernel estimator (see the discussion below).

To go further, for each sample size, we evaluate the mean squared error of the estimators, in other words
$$\mathrm{MSE}(\hat f)=\int\left(\hat{f}(x,y)-f(x,y)\right)^2dy,$$
where $\hat f$ is either the kernel rule or the projection estimate. In Appendix~\ref{simustables}, we give approximations of the MSE based on
$N=100$ samples for different values of $\eta$.

\begin{figure}[h]\begin{center}
 \includegraphics[scale=0.5]{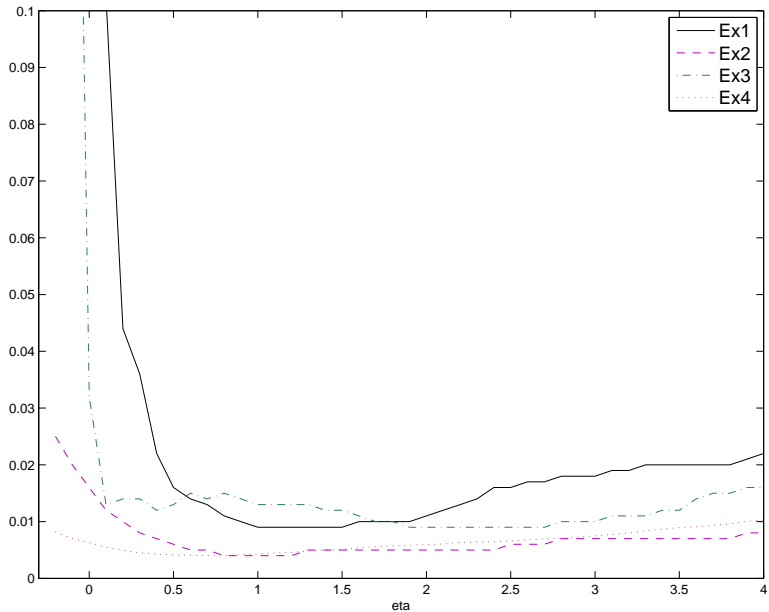}
  \end{center}
 \caption{MSE($\hat f$) for $n=500$, Example 1 ($x=0.5$), Example 2 ($x=0.5$), Example 3 ($x=0$), Example 4 ($x=0$)}
 \label{riskseloneta}
\end{figure}

Now, let us comment our results from the point of view of tuning, namely we try to answer the question: how to choose the parameter $\eta$? We first focus on kernel rules. Tables of Appendix~\ref{simustables} show that, often, the optimal value is $\eta=1$. More precisely, it is always the case for Examples 1 and 2. For Examples 3 and 4, when $\eta=1$ is not the optimal value, taking $\eta=1$ does not deteriorate the risk too much. So, for kernel rules, the choice $\eta=1$ is recommended even if larger values can be convenient in some situations.
To shed more lights on these numerical results,  in Figure~\ref{riskseloneta}, we draw the MSE for the kernel rule  in function of the parameter $\eta$.
We observe that the shape of the curve  is the same whatever the example. If $\eta$ is too small the risk blows up, which shows that the assumption $\eta>0$ in theoretical results is unavoidable at least asymptotically. Furthermore, we observe that if $\eta$ is too large, then the estimate oversmooths and the risk increases but without explosion for $\eta$ not too far from the minimizer.  Similar phenomena have already been observed for wavelet thresholding rules for density estimation (see Section 2.2 of \cite{densitywavelet}).  Tuning kernel rules is then achieved.

We now deal with projection rules. Unfortunately, the plateau phenomenon of Figure~\ref{riskseloneta} does not happen for projection estimators. In this case, the optimal value for $\eta$ seems to change according to the example. Tuning this procedure is not so obvious. Note that performances of kernel and projections rules are hardly comparable since they are respectively based on a Gaussian kernel function and piecewise constant functions.

For kernel rules, we study the influence of the knowledge of $\fx$. Tables \ref{table1.simul} and \ref{table3.simul} show that when $\fx$ is known results are a bit better as expected, but the difference is not very significant. Since projections rules are less sensitive to the estimate $\hfx$, we only show results with $\fx$ unknown.  Finally, to study the dependence of estimation with respect to $x$, we focus on Tables~\ref{table2.simul} and \ref{table2} that  show that in Example 3 estimation is better at $x=0$ and $x=1$ than at $x=0.36$. This was expected since the density design is smaller at $x=0.36$ and this confirms the role of $\delta$ in the rate of convergence of both estimators (see Theorems~\ref{theo2} and \ref{oracleproba}).
Similar conclusions can be drawn for Example 4.  Finally, we wish to mention that the ratio between the risk of our procedures and the oracle risk (the upper bounds of Theorems \ref{theo1},  \ref{theo2} and \ref{oracleproba}) remains bounded with respect to $n$, which corroborates our theoretical results.


\section{Proofs} \label{sec:proofs}
In this section, after giving intermediate technical results, we prove the results of our paper. Most of the time,  as explained in introduction, we only consider the case $d_1=d_2=1$. We use notations that we have previously defined. The classical Euclidian norm is denoted $\|\cdot\|$. 
Except if the context is ambiguous, from now on, the $\|\cdot\|_\infty$-norm shall denote the supremum either on $\R$, on $V_n(x)$ or on $V_n(x)\times\R$. We shall also use for any function $g$
$$\|g\|_{\infty,2}:=\sup_{t\in V_n(x)}\|g\|_{t,2}.$$
This section is divided into two parts: Section~\ref{sec:proofskernel} (respectively Section~\ref{sec:proofsprojections}) is devoted to the proofs of the results for the kernel rules (respectively for the projection rules). We first prove in Section \ref{sec:borneinforacle} the lower bound stated in Proposition \ref{borneinforacle}. Main results for kernel rules, namely Theorems \ref{theo1} and \ref{theo2} are proved in Section~\ref{sec:proof12}. They depend on several intermediate results that are proved in Sections~\ref{sec:proof3}--\ref{sec:proof5} (see the sketch of proofs in Section~\ref{sec:proof12}). Theorem~\ref{propvitesse} that derives rates for kernel rules is proved in Section~\ref{sec:prooftheo5}. For projection rules, the main theorem, namely Theorem~\ref{oracleproba}, is proved in Section~\ref{sec:prooforacleproba}. It is based on intermediate results shown in Sections~\ref{preuveconcentration}--\ref{preuvecoarsebound}. Finally, Theorem~\ref{vitesse} that derives rates for projection rules is proved in Section~\ref{SecPreuvebiais}. 
As usual in nonparametric statistics, our results are based on sharp concentration inequalities that are stated in Lemmas~\ref{bernstein}, \ref{tala2} and \ref{lemmafx}. These lemmas and other technical results stated in Lemmas~\ref{sup1} and \ref{majbias} are proved in Appendix~\ref{preuvestechniques}. 
\begin{lemma}\label{bernstein}
[Bernstein Inequality]
 Let $(U_i)$ be a sequence of i.i.d. variables uniformly bounded by a positive constant $c$ and such that
$\E U_1^2\leq v. $ Then
$$ \P\left(\left|\frac{1}{n}\sum_{i=1}^nU_i-\E[U_i]\right|\geq \varepsilon\right)\leq
2\exp\left(-\min\left(\frac{n\varepsilon^2}{4v},\frac{n\varepsilon}{4c}\right)\right)$$
\end{lemma}
Note that Lemma~\ref{bernstein} is a simple consequence of \citet{BM98}, p.366.

\begin{lemma}\label{tala2}
[Talagrand Inequality]
Let $U_1,\dots,U_n$ be i.i.d. random variables and  $\nu_n(a)=\frac 1n\sum_{i=1}^n [\tau_a(U_i)-{\mathbb E}(\tau_a(U_i))]$
for $a$ belonging to ${\mathcal A}$ a countable subset of functions. For any $\zeta>0$,
\begin{equation*}
 \P(\sup_{a\in {\mathcal A}}|\nu_n(a)|\geq (1+2\zeta)H)\leq 2\max\left(\exp\left(-\frac{\zeta^2}{6}\frac{nH^2}{v}\right),
\exp\left(-\frac{\min(\zeta,1)\zeta}{21}\frac{n H}{M}\right)\right)
\end{equation*}
with
$$\sup_{a\in   {\mathcal A}}\sup_u|\tau_a(u)|\leq M, \;\;\;\;
\mathbb{E}\Big[\sup_{a\in  {\mathcal A}}|\nu_n(a)|\Big]\leq H,
\; \;\;\;\sup_{a\in   {\mathcal A}}{\rm Var}(\tau_a(U_1)) \leq v.$$
\end{lemma}

Let $\rho>1$ and consider the event
$$\Lambda_\rho=\{
 \rho^{-1}\delta\leq \hat\delta\leq \rho\delta\}\cap\{\rho^{-2}\|\fx\|_\infty\leq \|\hfx\|_\infty  \leq \rho^{2}\|\fx\|_\infty
\}.$$
We have the following lemma.
\begin{lemma}\label{lemmafx}
Condition (\ref{Cfx})  implies that  $$\P(\Lambda_\rho^c)\leq B_1e^{-(\log n)^{3/2}}$$ with some positive constant $B_1$ that depends on $\fx$ and $\rho$.
\end{lemma}

\begin{lemma}\label{sup1}
For any integrable functions $f_1$ and $f_2$,
if the support of $u\mapsto f_2(u,y)$ is included in $[-2A/k_n,2A/k_n]^{d_1}$ for all $y$, then
we have
$$\|f_1*f_2\|_{x,2}\leq \sup_{t\in V_n(x)} \|f_1\|_{t,2}\times \|f_2\|_1,$$
\end{lemma}
\begin{lemma}\label{majbias} We use notations of Definition \ref{defSm}.
Let $m=(m_1,m_2)$ be fixed. For any function $\tau$, the projection $\mathcal{K}_{m}(\tau)$ of $\tau$ on $S_{m}$ verifies
$$\|\mathcal{K}_{m}(\tau)\|_{x,2}\leq (r+1) \|\fx\|_{\infty} \delta^{-1} \sup_{t\in V_{n}(x)}\|\tau\|_{t,2}.$$
\end{lemma}
\subsection{Proofs for the kernel estimator}\label{sec:proofskernel}
\subsubsection{Proof of Proposition~\ref{borneinforacle}}\label{sec:borneinforacle}
We just need to control:
\begin{eqnarray*}
\int\var(\hat{f}_{h}(x,y))dy&=&\frac{1}{n}\int \var\left(\left[\fx(X_1)\right]^{-1}K_{h}(x-X_1,y-Y_1)\right)dy\\
&=&\frac{1}{n}\int\left(  \E\left[\left[\fx(X_1)\right]^{-2}K_{h}^2(x-X_1,y-Y_1)\right]- \left(\E\left[\left[\fx(X_1)\right]^{-1}K_{h}(x-X_1,y-Y_1)\right]\right)^2\right)dy.
\end{eqnarray*}
First, by using Lemma~\ref{sup1} and (\ref{sup>2}),
\begin{eqnarray*}
\int  \left(\E\left[\left[\fx(X_1)\right]^{-1}K_{h}(x-X_1,y-Y_1)\right]\right)^2dy&=&\int (K_h*f)^2(x,y)dy\\
&\leq& \|K_h\|_1^2\times\sup_{t\in V_n(x)} \|f\|_{t,2}^2
\leq\|K\|_1^2 \|f\|_\infty.
\end{eqnarray*}
Furthermore,
\begin{eqnarray*}
\int  \E\left[\left[\fx(X_1)\right]^{-2}K_{h}^2(x-X_1,y-Y_1)\right]dy&=&\iiint K_h^2(x-u,y-v)f(u,v)[\fx(u)]^{-1}dudvdy\\
&=&\int (K_{h_1}^{(1)})^2(x-u)[\fx(u)]^{-1}du \times \frac{\|K^{(2)}\|_2^2}{h_2}\\
&=&\frac{\|K^{(2)}\|_2^2}{h_1h_2}\times \int\frac{[K^{(1)}(s)]^2}{\fx(x-sh_1)}ds.
\end{eqnarray*}
Now
assume that $\fx$ is positive and continuous on a neighborhood of $x$. Since $\max\mathcal{H}_n\to 0$ when $n\to +\infty$, then $h_1 \to 0$. Then we have
\begin{eqnarray*}
\left|\int\frac{\fx(x)[K^{(1)}(s)]^2}{\fx(x-sh_1)}ds -\|K^{(1)}\|_2^2\right|&\leq&\int[K^{(1)}(s)]^2 \left|\frac{\fx(x)}{\fx(x-sh_1)}-1\right|ds\\
& \le&  \max_{|v|\le Ah_1}\left| \frac{\fx(x)}{\fx(x+v)}-1 \right|\int[K^{(1)}(s)]^2 ds =o(1).
\end{eqnarray*}

\subsubsection{Proof of Theorems~\ref{theo1} and \ref{theo2}}\label{sec:proof12} 
We introduce $$g(x,y)=\frac{f_{X,Y}(x,y)}{\hfx(x)}=\frac{\fx(x)}{\hfx(x)}f(x,y).$$
%
We consider the set $\Gamma=\Gamma_1\cap\Gamma_2$ where
$$\Gamma_1=\left\{ \forall h,h'\in\mathcal{H}_n: \left\|K_{h}*\hat f_{h'}-K_{h}*K_{h'}*g\right\|_{x,2}\le \frac{\chi_1}{\sqrt{\hat\delta nh'_1h'_2}}\quad \right\},$$
$$\Gamma_2=\left\{ \forall h'\in\mathcal{H}_n: \left\|\hat f_{h'}-K_{h'}*g \right\|_{x,2}\le \frac{\chi_2}{\sqrt{\hat\delta nh'_1h'_2}}\quad \right\}$$
and $$\chi_1=(1+\eta)\|K\|_1\|K\|_2,\quad \chi_2= (1+\eta)\|K\|_2.$$
We shall use following propositions that deal with the general case when $\fx$ is estimated by $\hfx$. When $\fx$ is known, it can easily be checked that these propositions also hold with $g$ replaced by $f$ and $\hat\delta$ by $\delta$.
We also use the set $\Lambda_{\rho}$ studied in Lemma~\ref{lemmafx} with $ \rho=(1+\eta/2)^2$. 

Let us give a sketch of the proof. The main steps for proving  Theorems~\ref{theo1} and \ref{theo2} are the following. We first prove an oracle inequality for the function $g$ on the set $\Gamma$ (Proposition~\ref{prop1}). Then, in Proposition~\ref{prop2}, we prove that the event $\Gamma$ occurs with large probability by using Lemma~\ref{tala2}. Finally, Proposition~\ref{prop3} studies the impact of replacing $g$ by $f$. Proposition~\ref{prop5} gives a polynomial control in $n$ of our estimate that is enough to control its risk on $\Gamma^c$ by using Proposition~\ref{prop2} and Lemma~\ref{lemmafx}. 
 
\begin{proposition}\label{prop1}
On the set $\Gamma$, we have the following result.
\begin{equation*}
\|\hat f-g\|_{x,2}\le \inf_{h\in\mathcal{H}_n}\left\{C_1\|K_h*g-g\|_{\infty,2}+C_2\frac{1}{\sqrt{\hat{\delta}nh_1h_2}}\right\},
\end{equation*}
where $C_1=1+2\|K\|_1$  and $C_2=(1+\eta)\|K\|_2(3+2\|K\|_1)$.
\end{proposition}

\begin{proposition}\label{prop2}
 Under $(H_1)$, $(H_3)$ and $(CK)$, we have:
\begin{equation*}
\P\left( \Gamma^c\cap \Lambda_\rho\right)\le C\exp\{-(\log n)^{5/4}\}
\end{equation*}
where $C$ depends on $K$, $\eta$ and $\|f\|_\infty$.
\end{proposition}

\begin{proposition}\label{prop3} Assume that $(H_1)$, $(H_2)$ and $(CK)$ are satisfied. On $\Lambda_{\rho}$:
\begin{align*}
\|K_h*g-g\|_{\infty,2}\le& \|K_h*f-f\|_{\infty,2} +C\delta^{-1}\|\hfx-\fx\|_{\infty},\\
\|g-f\|_{x,2}\le& C\delta^{-1}\|\hfx-\fx\|_{\infty},
\end{align*}
where $C$ depends on $\eta$, $K$, and  $\|f\|_\infty$.
\end{proposition}

\begin{proposition}\label{prop5}
Assume that  $(CK)$ is satisfied. For any $h\in\mathcal{H}_n$,
$$\|\hat f_h\|_{x,2}\leq \|K^{(1)}\|_\infty\|K^{(2)}\|_2(\log n)^{-3} n^{3/2}.$$
\end{proposition}

The first part of Theorem~\ref{theo1} can be deduced from Propositions~\ref{prop1} and \ref{prop2}. Note that in the case of Theorem~\ref{theo1}, since $\fx$ is known, $g=f$ and $\P(\Lambda_\rho)=1$. The second part of Theorem~\ref{theo1} is a consequence of Proposition~\ref{prop5}, (\ref{sup>2}) and (\ref{traj1}). Since
$$\|\hat{f}-f\|_{x,2}\leq \|\hat f-g\|_{x,2}+\|g-f\|_{x,2}$$
and
$$\Gamma\cap \Lambda_\rho=(\Gamma\cup \Lambda_\rho^c)\cap \Lambda_\rho,$$
the first part of Theorem~\ref{theo2} is a consequence of Propositions~\ref{prop1}, \ref{prop2} and \ref{prop3} combined with Lemma~\ref{lemmafx}.
The second part of Theorem~\ref{theo2} is a consequence of Proposition~\ref{prop5}, (\ref{sup>2})  and (\ref{traj2}).
\subsubsection{Proof of Proposition~\ref{prop1}}\label{sec:proof3}
We apply the GLM as explained in Section~\ref{explicationGLM} with $\hat f_h$ given in (\ref{estim1}) for estimating $g$,
$\mathcal{M}_{n}=\mathcal{H}_n$, $\|.\|=\|.\|_{x,2}$,
$\sigma(h)=\chi/\sqrt{\hat\delta nh_1h_2}$,
 and the operator $\mathcal{K}_{h}$ is the convolution product with $K_{h}.$ Note that (\ref{commute1}), (\ref{commute2}) and (\ref{E-S}) are satisfied but not (\ref{norme-op}).
 But we have:
 $$B(h)=\sup_{h'\in\mathcal{H}_n}\|\mathcal{K}_{h'}(g)-(\mathcal{K}_{h'}\circ \mathcal{K}_h)(g)\|_{x,2}\leq \|K\|_{1}\sup_{t\in V_n(x)}\|g- \mathcal{K}_h(g)\|_{t,2},$$
 using Lemma~\ref{sup1} and the equality $\|K_{h'}\|_{1}=\|K\|_{1}$.
Let us fix $h\in \mathcal{H}_n$. We obtain Inequality \eqref{or1} in our case:
$$\|\hat f-g\|_{x,2}\leq 2B(h)+2\sigma(h)+\|\hat f_h-\mathcal{K}_h(g)\|_{x,2}+\|g-\mathcal{K}_h(g)\|_{x,2}+2\xi(h)$$
with
 $$\xi(h)=\sup_{h'\in\mathcal{H}_n}\left\{\|(\hat f_{h'}-\mathcal{K}_{h'}(g))-(\mathcal{K}_{h'}(\hat f_h)-(\mathcal{K}_{h'}\circ \mathcal{K}_h)(g))\|_{x,2}-\sigma(h')\right\}_+.$$
But, on $\Gamma$, $\forall h,h'\in\mathcal{H}_n$, $ \|\hat f_{h'}- \mathcal{K}_{h'}(g)\|_{x,2}\le \chi_2/\sqrt{\hat\delta nh'_1h'_2}$
and $\|\mathcal{K}_{h'}(\hat f_h)-(\mathcal{K}_{h'}\circ \mathcal{K}_h)(g)\|_{x,2}\le \chi_1/\sqrt{\hat\delta nh'_1h'_2}$,
so that $\xi(h)=0$.
Then, on $\Gamma$,
\begin{eqnarray*}
\|\hat f-g\|_{x,2}&\leq &2B(h)+2\sigma(h)+\frac{\chi_2}{\sqrt{\hat\delta nh_1h_2}}+\|g-\mathcal{K}_h(g)\|_{x,2}\\
&\leq& (2 \|K\|_{1}+1)\sup_{t\in V_n(x)}\|g- \mathcal{K}_h(g)\|_{t,2}+\frac{2\chi+\chi_2}{\sqrt{\hat\delta nh_1h_2}}
\end{eqnarray*}
with $2\chi+\chi_2=2\chi_{1}+3\chi_{2}=(1+\eta)(2\|K\|_1+3)\|K\|_2.$

\subsubsection{Proof of Proposition~\ref{prop2} }
We respectively denote $\tilde{\P}$ and $\tilde{\E}$  the probability distribution and the expectation associated with $(X_1,Y_1),\ldots,(X_n,Y_n)$.
Thus
\begin{eqnarray*}
\Gamma_1=\left\{ \forall h,h'\in\mathcal{H}_n: \left\|\hat{f}_{h,h'}-\tilde\E\left[\hat{f}_{h,h'} \right]\right\|_{x,2}\le \frac{\chi_1}{\sqrt{\hat\delta nh'_1h'_2}} \right\},\\
\Gamma_2=\left\{ \forall h'\in\mathcal{H}_n: \left\|\hat{f}_{h'}-\tilde\E\left[\hat{f}_{h'} \right]\right\|_{x,2}\le \frac{\chi_2}{\sqrt{\hat\delta nh'_1h'_2}} \right\}.
\end{eqnarray*}

To prove Proposition ~\ref{prop2}, we study $\Gamma_1^c\cap \Lambda_\rho$ and $\Gamma_2^c\cap \Lambda_\rho$. So first, let assume we are on the event $\Lambda_\rho$. Note that on $\Lambda_\rho$, we have $\hat{\delta}^{-1}\le\rho\delta^{-1}$ and for all $u\in V_n(x)$, $|g(u,v)|\le f(u,v)\rho$ (see the proof of Lemma~\ref{lemmafx}).
We denote for any $x$, $y$, $u$ and $v$,
$$w(x,y,u,v)=[\hfx(u)]^{-1}(K_h*K_{h'})(x-u,y-v).$$
We  can then write:
$$\hat f_{h,h'}(x,y)=\frac{1}{n}\sum_{i=1}^n w(x,y,X_i,Y_i)$$
and with ${\mathcal B}$ the unit ball in $\L_2(\R)$ endowed with the classical norm and ${\mathcal A}$ a dense countable subset of ${\mathcal B}$,
\begin{eqnarray*}
 \left\| \hat{f}_{h,h'}-\tilde\E\left[\hat{f}_{h,h'}\right]\right\|_{x,2}&=&\sup_{a\in {\mathcal B}}\int a(y)\left(\hat f_{h,h'}(x,y)-\tilde\E[\hat f_{h,h'}(x,y)]\right)dy\\
&=&\sup_{a\in {\mathcal A}}\int a(y)\left(\hat f_{h,h'}(x,y)-\tilde\E[\hat f_{h,h'}(x,y)]\right)dy\\
&=&\sup_{a\in {\mathcal A}}\frac{1}{n}\sum_{i=1}^n\int a(y)\left[w(x,y,X_i,Y_i)-\tilde\E(w(x,y,X_i,Y_i))\right]dy.
\end{eqnarray*}
Hence, one will apply the inequality of Lemma \ref{tala2}
with $\tau_{a,x}(X_i,Y_i)=\int a(y)w(x,y,X_i,Y_i)dy$.
First, we have:
\begin{eqnarray*}
\left(\tilde\E\left[ \left\| \hat{f}_{h,h'}-\tilde\E\left[\hat{f}_{h,h'}\right]\right\|_{x,2}\right]\right)^2&\leq&\tilde\E\left[ \left\| \hat{f}_{h,h'}-\tilde\E\left[\hat{f}_{h,h'}\right]\right\|_{x,2}^2\right]\\
&=&\tilde\E\left[ \int\left(\hat f_{h,h'}(x,y)-\tilde\E[\hat f_{h,h'}(x,y)]\right)^2 dy\right]\\
&=&\int \var(\hat f_{h,h'}(x,y))dy\\
&=&\frac{1}{n}\int \var\left([\hfx(X_1)]^{-1}(K_h*K_{h'})(x-X_1,y-Y_1)\right)dy\\
&\leq&\frac{1}{n}\int \tilde\E\left([\hfx(X_1)]^{-2}(K_h*K_{h'})^2(x-X_1,y-Y_1)\right)dy\\
&\leq&\frac{1}{\hat\delta n}\iiint (K_h*K_{h'})^2(x-u,y-v)g(u,v)dudvdy.
\end{eqnarray*}
But we have
\begin{eqnarray*}
(K_h*K_{h'})^2(x-u,y-v)&=&\left(\iint K_{h'}(x-u-s,y-v-t)K_{h}(s,t)dsdt\right)^2\\
&\leq&\iint K_{h'}^2(x-u-s,y-v-t)|K_{h}(s,t)|dsdt\times \|K\|_1.
\end{eqnarray*}
 Therefore, since for any $u$, $\int f(u,v)dv=1$ and $K(x,y)=K^{(1)}(x)K^{(2)}(y)$,
 \begin{eqnarray*}
 &&\left(\tilde\E\left[ \left\| \hat{f}_{h,h'}-\tilde\E\left[\hat{f}_{h,h'}\right]\right\|_{x,2}\right]\right)^2
 \leq\frac{\|K\|_1}{\hat\delta n}\iint  |K_{h}(s,t)|\left(\iiint  K_{h'}^2(x-u-s,y-v-t)g(u,v)dudvdy\right)dsdt\\
&&=\frac{\|K\|_1}{\hat\delta n}\iint  |K_{h}(s,t)|\left(\int\left(\int\left(\int (K^{(1)}_{h_1'})^2(x-u-s)(K^{(2)}_{h_2'})^2(y-v-t)dy\right) g(u,v)dv\right)du\right)dsdt\\
&&\leq\frac{\|K\|_1\|K^{(2)}\|_2^2\rho}{\hat\delta nh_2'}\iint  |K_{h}(s,t)|\left(\int   (K^{(1)}_{h_1'})^2(x-u-s) du\right)dsdt\\
&&=\frac{\|K\|_1^2\|K^{(1)}\|_2^2\|K^{(2)}\|_2^2\rho}{\hat\delta nh_1'h_2'}=\frac{\|K\|_1^2\|K\|_2^2\rho}{\hat\delta nh_1'h_2'}.
\end{eqnarray*}
Consequently, we obtain $\tilde\E\left[ \left\| \hat{f}_{h,h'}-\tilde\E\left[\hat{f}_{h,h'}\right]\right\|_{x,2}\right]\leq H$, with
\begin{equation}
H= \frac{\|K\|_1\|K\|_2\rho^{1/2}}{\sqrt{\hat\delta nh_1'h_2'}}.
\end{equation}
Now, let us deal with $v$ which is an upper bound of $\sup_{a\in {\mathcal A}}\var\left(\tau_{a,x}(X_1,Y_1) \right)$.
\begin{eqnarray*}
\sup_{a\in {\mathcal A}}\var\left(\tau_{a,x}(X_1,Y_1) \right)
&\leq&\sup_{a\in {\mathcal A}}\tilde\E\left[ \left(\int a(y)w(x,y,X_1,Y_1)dy\right)^2\right]\\
&\leq&\sup_{a\in {\mathcal A}}\tilde\E\left[\int|w(x,y,X_1,Y_1)|dy\int a^2(y)|w(x,y,X_1,Y_1)|dy\right]\\
&\leq&\sup_{u,v}\int|w(x,y,u,v)|dy\, \sup_{y}\tilde\E[|w(x,y,X_1,Y_1)|].
\end{eqnarray*}
Now,
\begin{eqnarray*}
\sup_{u,v}\int|w(x,y,u,v)|dy&=&\sup_{u,v}\int\left|[\hfx(u)]^{-1}(K_h*K_{h'})(x-u,y-v)\right|dy\\
&\leq&\frac{1}{\hat\delta }\sup_{u,v}\int\left|\iint K^{(1)}_{h_1'}(x-u-s) K^{(2)}_{h_2'}(y-v-t)K_{h}(s,t)dsdt\right|dy\\
&\leq& \frac{1}{\hat\delta }\sup_{u,v}\iint |K_{h}(s,t)|\left( \int |K^{(1)}_{h_1'}(x-u-s)| | K^{(2)}_{h_2'}(y-v-t) |dy\right) dsdt\\
&\leq&\frac{\|K\|_1\|K^{(2)}\|_1\|K^{(1)}\|_\infty}{\hat\delta h_1'}
\end{eqnarray*}
and
\begin{eqnarray*}
 \sup_{y}\tilde\E[|w(x,y,X_1,Y_1)|]&=&\sup_y\iint |w(x,y,u,v)|f_{X,Y}(u,v)dudv\\
 &=&\sup_y\iint |(K_h*K_{h'})(x-u,y-v)|g(u,v)dudv\\
 &\leq&\|g\|_\infty\sup_y\iint\left(\iint |K_h(x-u-s,y-v-t)||K_{h'}(s,t)|dsdt\right)dudv\\
 &\leq&\|g\|_\infty\|K\|_1^2\leq \|f\|_{\infty}\rho\|K\|_1^2
\end{eqnarray*}
since on $\Lambda_\rho$,  $\|g\|_{\infty}\le \rho\|f\|_\infty$ and where  $\|g\|_{\infty}=\sup_{(t,v)\in V_n(x)\times\R} |g(t,v)|$.
Thus, we set
\begin{equation}
v=\frac{\|K\|_1^3\|K^{(2)}\|_1\|K^{(1)}\|_\infty\rho\|f\|_\infty}{\hat\delta h_1'}.
\end{equation}
Finally, we deal with $M$ which has to be an upper bound of  $\sup_{a\in {\mathcal A}}\sup_u\sup_v\left|\int a(y)w(x,y,u,v)dy\right|$
\begin{eqnarray*}
\sup_{a\in {\mathcal A}}\sup_u\sup_v\left|\int a(y)w(x,y,u,v)dy\right|
&=&\sup_{u,v}\|w(x,.,u,v)\|_2\\
&\leq&\frac{1}{\hat\delta}\sup_{u,v}\left(\int (K_h*K_{h'})^2(x-u,y-v)dy\right)^{1/2}.
\end{eqnarray*}
We have:
\begin{eqnarray*}
\int (K_h*K_{h'})^2(x-u,y-v)dy&=&\int\left(\iint K_{h'}(x-u-s,y-v-t)K_{h}(s,t)dsdt\right)^2dy\\
&\leq&\|K\|_1\iint|K_{h}(s,t)|\left(\int K_{h'}^2(x-u-s,y-v-t)dy\right)dsdt\\
&\leq&\frac{\|K\|_1^2\|K^{(1)}\|_\infty^2\|K^{(2)}\|_2^2}{h_1'^2h_2'}.
\end{eqnarray*}
Therefore, we can set
\begin{equation}
M=\frac{\|K\|_1\|K^{(1)}\|_\infty\|K^{(2)}\|_2}{\hat\delta h_1'\sqrt{h_2'}}.
\end{equation}
So, since $\rho=(1+\eta/2)^{2}$, Lemma \ref{tala2} implies that for any $\zeta>0$,
\begin{align*}&\tilde\P\left(\left\| \hat{f}_{h,h'}-\tilde\E\left[\hat{f}_{h,h'}\right]\right\|_{x,2}\geq  (1+2\zeta)\frac{ \|K\|_1\|K\|_2(1+\eta/2)}{\sqrt{\hat{\delta}nh_1'h_2'}}\right)\\&\leq 2\max\left(\exp\left\{ -\frac{\zeta^2C_1(K,\|f\|_\infty)}{h_2'} \right\} ,  \exp\left\{-\zeta\min(1,\zeta)C_2(K,\eta)\sqrt{nh_1'\hat\delta} \right\} \right),
\end{align*}
where $C_1(K,\|f\|_\infty)$ and $C_2(K,\eta)$ are positive constants that depend on $K$ and $\|f\|_\infty$ and $K$ and $\eta$ respectively.
Similarly we have for any $\zeta>0$,
\begin{align*}&\tilde\P\left(\left\| \hat{f}_{h'}-\tilde\E\left[\hat{f}_{h'}\right]\right\|_{x,2}\geq  (1+2\zeta)\frac{ \|K\|_2(1+\eta/2)}{\sqrt{\hat{\delta}nh_1'h_2'}}\right)\\&\leq 2\max\left(\exp\left\{ -\frac{\zeta^2C_3(K,\eta,\|f\|_\infty)}{h_2'} \right\} ,  \exp\left\{-\zeta\min(1,\zeta)C_4(K,\eta)\sqrt{nh_1'\hat\delta} \right\} \right),
\end{align*}
where $C_3(K,\eta,\|f\|_\infty)$ and $C_4(K,\eta)$ are positive constants that depend on  $K$, $\eta$ and $\|f\|_\infty$ and $K$ and $\eta$ respectively.
Let $\zeta=\eta/(4+2\eta)$ so that $(1+2\zeta)(1+\eta/2)=(1+\eta)$.
For $(h'_1,h'_2)\in\mathcal{H}_n$, $\frac{(\log n)^3}{\rho n}\leq\frac{(\log n)^3}{\rho \delta n}\leq\frac{(\log n)^3}{\hat{\delta}n}\leq h'_1<1$ and $\frac{1}{n} \leq h'_2<\frac{1}{(\log n)^2-1}$. So, $-\sqrt{nh_1'\hat\delta}\leq -(\log n)^{3/2}$ and $-\frac{1}{h'_2}<-(\log n)^2+1$. Therefore, on $\Lambda_\rho$,
\begin{eqnarray}
&&\sum_{h,h'\in\mathcal{H}_n}\tilde\P\left(\left\| \hat{f}_{h,h'}-\tilde\E\left[\hat{f}_{h,h'}\right]\right\|_{x,2}\geq  (1+\eta)\frac{ \|K\|_1\|K\|_2}{\sqrt{\hat{\delta}nh_1'h_2'}}\right)
\le \rho^2n^4e^{-C_5(K,\eta,\|f\|_\infty)(\log n)^{3/2}}\nonumber\\
&&\le  C_6(K,\eta,\|f\|_\infty)e^{-(\log n)^{5/4}},\label{proba2}
\end{eqnarray}
with $C_5(K,\eta,\|f\|_\infty)$ and $C_6(K,\eta,\|f\|_\infty)$ positive constants depending on $K$, $\eta$ and $\| f\|_\infty$. We have a similar result for $\sum_{h'\in\mathcal{H}_n}\tilde\P\left(\left\| \hat{f}_{h'}-\tilde\E\left[\hat{f}_{h'}\right]\right\|_{x,2}\geq  (1+\eta)\frac{ \|K\|_2}{\sqrt{\hat{\delta}nh_1'h_2'}}\right)$. Now to conclude, note that the right hand side of Inequality (\ref{proba2}) is not random. This allows us to obtain the result of the proposition.
\subsubsection{Proof of Proposition~\ref{prop3} }
We have the following decomposition
\begin{equation}\label{equ2}
K_h*g-g=K_h*g-K_h*f+K_h*f-f+f-g.
\end{equation}
Next,  on $\Lambda_\rho$,
\begin{align*}
\left|K_h*g(x,y)-K_h*f(x,y)\right|&=\left|\iint K_h(x-u,y-v)\left(g(u,v)-f(u,v)\right)dudv\right|\\
& = \left|\iint K_h(x-u,y-v)\frac{f(u,v)}{\hfx(u)}\left(\fx(u)-\hfx(u)\right)dudv\right|\\
& \le \sup_{t\in V_n(x)}\left|\fx(t)-\hfx(t)\right| \hat{\delta}^{-1} \iint |K_h(x-u,y-v)|f(u,v)dudv\\
&\le \sup_{t\in V_n(x)}\left|\fx(t)-\hfx(t)\right| \delta^{-1} \rho \iint |K_h(x-u,y-v)|f(u,v)dudv.
\end{align*}
Now by using (\ref{sup>2}), we have:
\begin{align*}
&\int\left(\iint |K_h(x-u,y-v)|f(u,v)dudv\right)^2dy\le \|K\|_1\iiint |K_h(x-u,y-v)|f^2(u,v)dudvdy\\
&\le  \|K\|_1\|K^{(2)}\|_1\iint |K^{(1)}_{h_1}(x-u)|f^2(u,v)dudv\le \|f\|_{\infty} \|K\|_1^2.
\end{align*}
Then we deduce that
\begin{equation}\label{equ1}
\|K_h*g-K_h*f\|_{\infty,2}\le C\delta^{-1}\sup_{t\in V_n(x)}\left|\fx(t)-\hfx(t)\right|,
\end{equation}
where $C$ depends on $\rho$, $\|f\|_{\infty}$ and $K$.
Moreover we have on $\Lambda_\rho$:
\begin{align*}\|g-f\|^2_{t,2}= &\int \frac{f^2(t,y)}{\hfx^2(t)}\left(\hfx(t)-\fx(t)\right)^2dy\le \|f\|_{\infty}\hat{\delta}^{-2}|\hfx(t)-\fx(t)|^2 \\
\le & C\delta^{-2}  |\fx(t)-\hfx(t)|^2,
\end{align*}
where $C$ depends on $\rho$ and $\|f\|_{\infty}$.
The last line, (\ref{equ2}) and (\ref{equ1}) allow us to conclude.
\subsubsection{Proof of Proposition~\ref{prop5}}\label{sec:proof5}
For any $h\in\mathcal{H}_n$, we have
$$\frac{1}{n\hat\delta h_1}\leq \frac{1}{(\log n)^3},\quad \frac{1}{h_2}\leq n.$$
Therefore,
\begin{eqnarray*}
\|\hat f_h\|_{x,2}^2&\leq&\int\left(\frac{1}{n}\sum_{i=1}^{n}\left|\hfx(X_i)\right|^{-1}\frac{1}{h_1}\left|K^{(1)}\left(\frac{x-X_i}{h_1}\right) \right|  \frac{1}{h_2}\left|K^{(2)}\left(\frac{y-Y_i}{h_2}\right) \right| \right)^2dy\\
&\leq&\int\left(\frac{1}{n}\sum_{i=1}^{n}\frac{1}{\hat\delta h_1}\|K^{(1)}\|_\infty  \frac{1}{h_2}\left|K^{(2)}\left(\frac{y-Y_i}{h_2}\right) \right| \right)^2dy\\
&\leq&n(\log n)^{-6}\|K^{(1)}\|_\infty^2 \sum_{i=1}^n\int \frac{1}{h_2^2}\left|K^{(2)}\left(\frac{y-Y_i}{h_2}\right) \right|^2dy\\
&\leq&n^3(\log n)^{-6}\|K^{(1)}\|_\infty^2 \|K^{(2)}\|_2^2,
\end{eqnarray*}
which proves the result.
\subsubsection{Proof of Theorem ~\ref{propvitesse}}\label{sec:prooftheo5}
We first assume that $d_1=d_2=1$. Using conditions $(BK_{\bold M})$, we then have:
\begin{eqnarray*}
(K_h*f)(x,y)-f(x,y) &=& \iint K(u,v)\left[f(x-uh_1,y-vh_2)-f(x,y)\right]dudv\\
&=&\iint K(u,v)\left[f(x-uh_1,y-vh_2)-f(x,y-vh_2)+f(x,y-vh_2)-f(x,y)\right]dudv\\
&= &\iint K(u,v)\left[\frac{(-uh_1)^{\lfloor \alpha_1\rfloor}}{\lfloor \alpha_1\rfloor!}\left(\frac{d^{\lfloor \alpha_1\rfloor}}{dx^{\lfloor \alpha_1\rfloor}}f(x+\tilde{u}h_1,y-vh_2)-\frac{d^{\lfloor \alpha_1\rfloor}}{dx^{\lfloor \alpha_1\rfloor}}f(x,y-vh_2)\right)\right]dudv\\
&+&\iint K(u,v)\left[\frac{(-vh_2)^{\lfloor \alpha_2\rfloor}}{\lfloor \alpha_2\rfloor!}\left(\frac{d^{\lfloor \alpha_2\rfloor}}{dy^{\lfloor \alpha_2\rfloor}}f(x,y+\tilde{v}h_2)-\frac{d^{\lfloor \alpha_2\rfloor}}{dy^{\lfloor \alpha_2\rfloor}}f(x,y)\right)\right]dudv\\
\end{eqnarray*}
where  $|\tilde{u}|\le |u|$ and $|\tilde{v}|\le |v|$.
If $f\in\mathcal H_{2}(\boldsymbol\alpha,\bold L)$, this implies that
\begin{equation*}
|(K_h*f)(x,y)-f(x,y)|\le C_1L_1h_1^{\alpha_1}+C_2L_2h_2^{\alpha_2},
\end{equation*}
where $C_1$ and $C_2$ depend on $\alpha_1$, $\alpha_2$ and $K$.
We can easily generalize this result to the case $d_1,d_2\geq 2$ and we obtain:
\begin{equation*}\label{eqpropbias}
|(K_h*f)(x,y)-f(x,y)|\le C\sum_{i=1}^d L_ih_i^{\alpha_i},
\end{equation*}
with a constant $C$ depending on $\boldsymbol{\alpha}$ and $K$.
Now taking
$$h_i=L_i^{-\frac{1}{\alpha_i}}\Delta_n^{-\frac{1}{\alpha_i}},\quad \Delta_n=\left(\prod_{i=1}^dL_i^{\frac{1}{\alpha_i}}\right)^{-\frac{\overline{\alpha}}{2\overline{\alpha}+1}}\left(\delta n\right)^{\frac{\overline{\alpha}}{2\overline{\alpha}+1}},$$
we obtain that
$$\frac{1}{\sqrt{\delta n\prod_{i=1}^d h_i}}=\Delta_n^{-1}$$
and
\begin{equation*}\label{eqbias}\sup_{t\in V_n(x)}\|K_h*f-f\|_{t,2}\le C(\delta n)^{-\frac{\bar\alpha}{2\bar\alpha+1}}\left(\prod_{i=1}^dL_i^{\frac{1}{\alpha_i}}\right)^{\frac{\bar\alpha}{2\bar\alpha+1}},
\end{equation*}
using $(H_4)$ and where $C$ is a positive constant that does not depend on $\delta$, $n$ and $\bold{L}$. By using Theorem~\ref{theo2}, this concludes the proof of Theorem ~\ref{propvitesse}.
\subsection{Proofs for the projection estimator}\label{sec:proofsprojections}
The structure of the proof of the main theorem, namely Theorem~3, is similar to the structure of the proofs for kernel rules. It is detailed along Section \ref{sec:prooforacleproba}.

\subsubsection{Proof of Theorem~\ref{oracleproba}}\label{sec:prooforacleproba}

First, let $$\Gamma=\{\forall m\in \mathcal{M}_n\quad\|\hat f_{m}- \mathcal{K}_m(f)\|_{x,2}\leq \sigma(m)/2\}.$$
To prove Theorem~\ref{oracleproba}, we follow the GLM, as explained in Section~\ref{explicationGLM}, with
$\|.\|=\|.\|_{x,2}$,
 and the operator $\mathcal{K}_{m}$ is the projection on $S_{m}$.
  In this case, using Lemma~\ref{majbias},
 $$B(m)=\sup_{m'\in\mathcal{M}_n}\|\mathcal{K}_{m'}(f)-(\mathcal{K}_{m'}\circ \mathcal{K}_m)(f)\|_{x,2}\leq (r+1)\|\fx\|_{\infty}\delta^{-1}\sup_{t\in V_n(x)}\|f- \mathcal{K}_m(f)\|_{t,2}.$$
Moreover for all $m,m'\in \mathcal{M}_{n}$,  $\mathcal{K}_{m'}\circ \mathcal{K}_m=\mathcal{K}_{m\wedge m'}$,
with $m \wedge m'=(\min(m_1,m'_1),\min(m_2,m'_2))$,
and  $\sigma(m \wedge m')\leq \sigma(m')$.
As already explained in Section~\ref{explicationGLM}, we introduce
$\tilde{\mathcal{K}}_m(\hat{f}_{m'})=\hat{f}_{m\wedge m'}$ and
$$\xi(m)=\sup_{m'\in\mathcal{M}_n}\left\{\|(\hat f_{m'}-\mathcal{K}_{m'}(f))-(\tilde{\mathcal{K}}_{m'}(\hat f_m)-(\mathcal{K}_{m'}\circ \mathcal{K}_m)(f))\|_{x,2}-\sigma(m')\right\}_+.$$
Let us fix $m\in \mathcal{M}_{n}$. We obtain inequality \eqref{or1} in our case:
$$\|\tilde f-f\|_{x,2}\leq 2B(m)+2\sigma(m)+\|\hat f_m-\mathcal{K}_m(f)\|_{x,2}+\|f-\mathcal{K}_m(f)\|_{x,2}+2\xi(m).$$
But, on $\Gamma$,
for all $m, m'$ in $\mathcal{M}_{n}$, $ \|\hat{f}_{m'}- \mathcal{K}_{m'}(f)\|_{x,2}\le \sigma(m')/2$ and
$ \|\hat{f}_{m\wedge m'}- \mathcal{K}_{m\wedge m'}(f)\|_{x,2}\le \sigma (m')/2$,
so that $\xi(m)=0$.
Then, on $\Gamma$,
\begin{eqnarray}\label{ineq}
\|\tilde f-f\|_{x,2}&\leq &2B(m)+2\sigma(m)+\frac{\sigma(m)}{2}+\|f-\mathcal{K}_m(f)\|_{x,2}\nonumber\\
&\leq& (2(r+1)\|\fx\|_{\infty}\delta^{-1} +1)\sup_{t\in V_n(x)}\|f-\mathcal{K}_m(f)\|_{t,2}+\frac52\sigma(m).
\end{eqnarray}
Now, let $\|.\|_n$ be the empirical norm defined by
\begin{equation*}
\|t\|_n=\left(\frac{1}{n}\sum_{i=1}^{n}
t^2(X_i)\right)^{1/2} \label{empiricalnorm}
\end{equation*}
and $l_{m_1}$  be the index  such that $x$ belongs to the interval $I_{l_{m_1}}$. For $\rho=(1+\eta)^{1/5}$, let
\begin{eqnarray*}
\displaystyle \Omega_\rho
&=&\left\{\forall\, m,\quad  \forall t\in {\rm Span}(\varphi_{l_{m_1},d}^m)_{0\leq d\leq r} \quad\|t\|_n^2\geq \rho^{-1}\int t^2(u)\fx(u)du\right\}.
\end{eqnarray*}
The heart of the proof of Theorem~\ref{oracleproba}
is the following concentration result:
\begin{proposition} \label{concentrationbis }
Assume that assumptions $(H_1)$, $(H_2)$, $(H_3)$ and $(CM)$ are satisfied.
 There exists $C>0$ only depending on $\eta, \phi_1, \phi_2, r$, $\|f\|_\infty$ and $\|\fx\|_\infty$ and  $\delta$ such that
$$\P\left( \Gamma^c\cap \Lambda_\rho\cap\Omega_\rho\right)\le C\exp\{-(\log n)^{5/4}\}.$$
\end{proposition}
Proposition~\ref{concentrationbis } and the following result show that the event $\Gamma$ occurs with large probability.
\begin{proposition}\label{omegac} Assume that assumptions $(H_2)$, $(H_3)$ and $(CM)$ are satisfied. Then,
$$\P(\Omega_{\rho}^{c}\cap\Lambda_\rho)\leq
C\exp\{-(\log n)^{5/4}\},$$
where
$C$ is a constant only depending on $\rho, \phi_1, r, \|\fx\|_\infty$ and $\delta.$
\end{proposition}
Then, using  Lemma \ref{lemmafx} and Propositions \ref{concentrationbis } and \ref{omegac},
\begin{eqnarray}\label{ineq2}
 \P(\Gamma^c)\leq\P((\Gamma\cap\Lambda_\rho\cap \Omega_\rho)^{c})
 = \P(\Gamma^c\cap\Lambda_\rho\cap \Omega_\rho)+\P(\Omega_\rho^c\cap\Lambda_\rho)+\P(\Lambda_\rho^c)
  \leq K e^{-\log^{5/4}(n)}
\end{eqnarray}
with $K$ depending on
$\eta, \phi_1, \phi_2, r$, $\|f\|_\infty$ and $\fx$.
Then, the first part of Theorem~\ref{oracleproba} is proved.
To deduce the second part, we use the following proposition.
\begin{proposition}\label{coarsebound} 
For all $m\in\mathcal{M}_n$,
 $$\|f-\hat f_m\|_{x,2}^2\leq
 2\|f\|_\infty +2 (1+\eta)^{4/5}{\hat{\delta}}^{-2}(r+1)\phi_1^{2}\phi_2D_{m_1}^{2}D_{m_2}^{2}.$$
\end{proposition}
Using assumption $(CM)$, it implies that
 $\|f-\hat f_{\hat{m}}\|_{x,2}^2\leq \tilde C_3^2n^4,$
 where $\tilde C_3$ depends on $\eta, r, \phi_1, \phi_2$ and $\|f\|_\infty$.
Then, by using (\ref{ineq}) which is true on $\Gamma\cap \Lambda_{\rho}$ we have
\begin{eqnarray*}
 \E\|\tilde f-f\|_{x,2}^q&=&\E\|\tilde f-f\|_{x,2}^q\1_{\Gamma\cap \Lambda_{\rho}}+ \E\|\tilde f-f\|_{x,2}^q\1_{(\Gamma\cap \Lambda_{\rho})^c}\\
 &\leq&\tilde C_4\left(\sup_{t\in V_n(x)}\|f-\mathcal{K}_m(f)\|_{t,2}+\sqrt{\frac{\|\fx\|_\infty}{\delta}}\sqrt{\frac{D_{m_1}D_{m_2}}{\delta n}}
\right)^q+ \tilde C_3^qn^{2q}\P((\Gamma\cap\Lambda_\rho\cap \Omega_\rho)^{c}),
 \end{eqnarray*}
where $\tilde C_4$ depends on $\eta, \phi_1, \phi_2, r, \|\fx\|_\infty$ and $\delta$.
Using (\ref{ineq2}), this concludes the proof of Theorem~\ref{oracleproba}.
\subsubsection{Proof of Proposition~\ref{concentrationbis }}\label{preuveconcentration}
First, we introduce some preliminary material.
For any matrix $M$, we denote
$$\|M\|_{2}=\sup_{x\neq 0}\frac{\|Mx\|}{\|x\|}, \qquad \|M\|_{F}=\left(\sum_{j,k} |M_{j,k}|^{2}\right)^{\frac{1}{2}} $$
the operator norm and the Frobenius norm.
We shall use that for any matrices $M$ and $N$,
$$\|M\|_{2}\leq\|M\|_{F},\quad
\|MN\|_{2}\leq\|M\|_{2}\|N\|_{2},\quad
\|MN\|_{F}\leq\|M\|_{2}\|N\|_{F}.
$$
Now we fix $m\in \mathcal{M}_n$. Then the index $l_{m_1}$ such that $x$ belongs to the interval $I_{l_{m_1}}$ is fixed. For the sake of simplicity, we denote it by $l$. Note that $I_{l}\subset V_{n}(x)$, since $2^{-m_1}\leq k_{n}^{-1}$. We set
$$F_{m_1}^{(l)}={\rm Span}(\varphi_{l,d}^{m})_{0\leq d\leq r}.$$
Moreover we denote
$$\hat{G}= \hat{G}_m^{(l)},\quad \hat{Z}=\hat{Z}_m^{(l)}, \quad \hat{A}=\hat{A}_m^{(l)},\quad\varphi_d=\varphi_{l,d}^m, \quad \psi_k=\psi_k^{m}.$$
The elements of $\hat A$  are denoted $(\hat a_{d,k})_{d,k}$  instead of $(\hat a_{(l_{m_1},d),k}^m)_{d,k}$.
We also introduce
$${G}=\E(\hat{G})=\left(\langle
\varphi_{d_1},\varphi_{d_2}
\rangle_{X}\right)
_{0\leq d_1, d_2\leq r} $$
and
$$Z=\E(\hat{Z})=\left(\iint
\varphi_{d}(u)\psi_k(y)
f(u,y)\fx(u)dudy\right) _{0\leq d\leq r, k \in K_m}.$$
By using Lemma~\ref{obvious}, the coefficients $(a_{j,k}^m)$ of $\mathcal{K}_m(f)$ in the basis  verify the matrix equation $GA=Z$ where the coefficients of the matrix $A$ are $A_{d,k}=a_{(l_{m_1},d),k}^m$ but are denoted $a_{d,k}$ for short.
We shall use the following algebra result. If $M$ is a symmetric matrix,
$$\min({\rm Sp}(M))=\min_u \frac{u^{*}Mu}{u^{*}u}.$$
Then
\begin{equation}\label{min1}
\min(\Sp(G))=\min_u \frac{u^{*}Gu}{u^{*}u}=\min_{t\in F_{m_1}^{(l)}}\frac{\int t^2(u)\fx(u)du}{\|t\|_2^2}\geq \delta
\end{equation}
and, in the same way,
$$\min(\Sp(\hat{G}))=\min \frac{u^{*}\hat{G}u}{u^{*}u}=\min_{t\in F_{m_1}^{(l)}
}\frac{\|t\|_{n}^2}{\|t\|_2^2},$$
so that
\begin{equation}\label{min3}
 \text{ on }\Omega_{\rho}\quad \min(\Sp(\hat{G}))\geq \rho^{-1}\delta.
\end{equation}

Now, let us begin the proof of Proposition~\ref{concentrationbis }.
Since
$$(\hat f_{m}- \mathcal{K}_m(f))(x,y)=\sum_{d=0}^r\sum_{k\in K_m}(\hat{a}_{d,k}- a_{d,k})\varphi_{d}(x)\psi_k(y)$$
we deduce
\begin{eqnarray*}
\|\hat f_{m}- \mathcal{K}_m(f)\|_{x,2}^{2}&=&\sum_{k}(\sum_{d}(\hat{a}_{d,k}-a_{d,k})\varphi_{d}(x))^{2}
\leq \sum_{d}\varphi_{d}^{2} (x) \sum_{k}\sum_{d}(\hat{a}_{d,k}-a_{d,k})^{2}\\
&\leq & \phi_1 D_{m_1}\|\hat{A}-A\|^{2}_{F}.
\end{eqnarray*}
On $\Lambda_{\rho}$, $\delta\geq\rho^{-1}\hat{\delta}$.
Then, using \eqref{min3}, on $\Omega_{\rho}\cap\Lambda_{\rho}$, $\min(\Sp(\hat{G}))\geq \rho^{-2}\hat{\delta}=(1+\eta)^{-2/5}\hat{\delta}$,
so   we are in the case where
$\hat{A}=\hat{G}^{-1}\hat{Z}$. From now on, we always assume that we are on $\Omega_{\rho}\cap\Lambda_{\rho}$.
We have:
\begin{eqnarray*}
\|\hat{A}-A\|_{F} &\leq &\|(\hat{G}^{-1}-G^{-1})Z\|_{F}+\|\hat{G}^{-1}(\hat{Z}-Z)\|_{F}\\
&\leq &\|\hat{G}^{-1}-G^{-1}\|_{2} \|Z\|_{F}+\|\hat{G}^{-1}\|_{2} \|\hat{Z}-Z\|_{F}.
\end{eqnarray*}
Since $\hat{G}$ is symmetric, $\|\hat{G}^{-1}\|_{2}$ is equal to the spectral radius of $\hat{G}^{-1}$.
And, using \eqref{min3}, its eigenvalues are positive, then
$$\|\hat{G}^{-1}\|_{2}
=(\min (\Sp(\hat{G})))^{-1}\leq \rho \delta^{-1}.$$
In the same way, using \eqref{min1},
$$\|G^{-1}\|_{2}
=(\min (\Sp(G)))^{-1}\leq  \delta^{-1}.$$
Then,
$$\|\hat{G}^{-1}-G^{-1}\|_{2}=
\|\hat{G}^{-1}(G-\hat{G})G^{-1}\|_{2}\leq \rho \delta^{-2}\|G-\hat{G}\|_{2}
\leq \rho \delta^{-2}\|G-\hat{G}\|_{F}.$$
Thus
\begin{eqnarray*}
\|\hat{A}-A\|_{F}
&\leq &\rho \delta^{-2}\|G-\hat{G}\|_{F} \|Z\|_{F}+ \rho {\delta}^{-1}\|\hat{Z}-Z\|_{F}.
\end{eqnarray*}
Moreover, since for any function $s$, $\sum_d \langle s,\varphi_d \rangle^2\leq \int_{I_l} s^2(u)du$, where $\langle,\rangle$ denotes the standard $\L_2$ dot product,
\begin{eqnarray*}
\|Z\|_{F}^{2}
&=& \sum_{d=0}^r\sum_{k\in K_m}\langle \int \varphi_{d}(u)f(u,.)\fx(u)du , \psi_k\rangle^2\\
&\leq &\sum_{d=0}^r\int \left(\int \varphi_{d}(u)f(u,y)\fx(u)du\right)^{2}dy\\
&\leq &\int \int_{I_l} f^2(u,y)\fx^2(u)dudy\\
&\leq & \|\fx\|_{\infty}^{2}\|f\|_{\infty} (4A2^{-m_1}).
\end{eqnarray*}
Finally (still on $\Omega_{\rho}\cap \Lambda_\rho$),
\begin{eqnarray*}
\|\hat f_{m}- \mathcal{K}_m(f)\|_{x,2}&\leq & C_3\|\hat{G}-G\|_{F} +\rho{\delta}^{-1}\sqrt{\phi_1D_{m_1}}  \|\hat{Z}-Z\|_{F}.
\end{eqnarray*}
Here $C_3=\|\fx\|_\infty\rho\delta^{-2}(r+1)\sqrt{\|f\|_{\infty} }$.
Thus, with $\P_\rho(\cdot)=\P(\cdot\cap\Lambda_\rho\cap\Omega_\rho)$, we can write:
$$\P_\rho\left(\|\hat f_{m}- \mathcal{K}_m(f)\|_{x,2}\geq \frac{\sigma(m)}{2}\right)\leq P_{1,m}+P_{2,m}$$
with
$$\begin{cases}
 \displaystyle   P_{1,m}=\P_\rho\left(\|\hat{Z}-Z\|_{F}\geq \frac{\sigma(m)}{2\rho^2{\delta}^{-1} \sqrt{\phi_1D_{m_1}}}\right)\\
  \displaystyle   P_{2,m}=\P_\rho\left(\|\hat{G}-G\|_{F}\geq \frac{\sigma(m)}{2\rho C_3}(\rho -1)\right).
  \end{cases}$$

\bigskip

\noindent$[1]$ \emph{Study of $P_{1,m}$}:
Let $\nu_{n}(t)=\frac1n \sum_{i=1}^{n}t(X_i,Y_i)-\E(t(X_i,Y_i))$ and
$$S_m^{(l)}=F_{m_1}^{(l)}\otimes H_{m_2}=
\left\{t, \quad t(x,y)=\sum_{d=0}^r\sum_{k\in K_m} {b}_{d,k} \varphi_d(x)\psi_k(y), \; b_{d,k}\in {\mathbb R}\right\}.$$
Then,
\begin{eqnarray*}
 \sup_{t\in S_m^{(l)}, \|t\|_{2}\leq 1}|\nu_n(t)|^{2}&=&
\sum_{d,k}\left|\nu_n(\varphi_d\otimes\psi_k)\right|^{2}\\
&=&
\sum_{d,k}\left| \frac1n \sum_{i=1}^{n}\varphi_{d}(X_i)\psi_k(Y_i)-\E(\varphi_{d}(X_i)\psi_k(Y_i))\right|^{2}\\
&=&\|\hat{Z}-Z\|_{F}^{2}.
\end{eqnarray*}

We are reduced to bound:
$$
\P_\rho\left(\sup_{ t\in S_m^{(l)}, \|t\|_2\leq 1}|\nu_n(t)|\geq\frac{\sigma(m)}{2\rho^2{\delta}^{-1} \sqrt{\phi_1D_{m_1}}}\right).$$
To deal with this term, we use Lemma \ref{tala2}. So, we consider ${\mathcal A}$ a dense subset of $\{ t\in S_m^{(l)}, \|t\|_2\leq 1\}$ and we compute $M,H$ and $v$.

 $\bullet$ First,  if $t=\sum_{d,k}b_{dk}\varphi_{d}\otimes\psi_k$ then
 $$|t(u,v)|^2=|\sum_{d,k} b_{dk}\varphi_d(u)\psi_k(v)|^2
 \leq \sum_{d,k} b_{dk}^2\sum_{d,k}|\varphi_d(x)\psi_k(v)|^2
 \leq \|t\|^2_{2}\phi_1D_{m_1}\phi_2D_{m_2}.
 $$
 Thus $\sup_{t\in   {\mathcal A}}\|t\|_{\infty}\leq \sqrt{\phi_1\phi_2D_{m_1}D_{m_2}}$
 and we can take $ M=\sqrt{\phi_1\phi_2D_{m_1}D_{m_2}}$.

 $\bullet$ Secondly, we recall
 \begin{eqnarray*}
 \sup_{t\in S_m^{(l)}, \|t\|_{2}\leq 1}|\nu_n(t)|^{2}=
\sum_{d,k}\left| \frac1n \sum_{i=1}^{n}\varphi_{d}(X_i)\psi_k(Y_i)-\E(\varphi_{d}(X_i)\psi_k(Y_i))\right|^{2}.
  \end{eqnarray*}
 Since the data are independent,
  \begin{align*}
 \text{Var}\left(\frac{1}{n}\sum_{i=1}^n\varphi_{d}(X_i)\psi_k(Y_i)\right)
=\frac{1}{n}\text{Var}(\varphi_{d}(X_1)\psi_k(Y_1)).
\end{align*}
We deduce:
\begin{eqnarray*}
\sum_k\E\left| \frac1n \sum_{i=1}^{n}\varphi_{d}(X_i)\psi_k(Y_i)-\E(\varphi_{d}(X_i)\psi_k(Y_i))\right|^{2}&\leq &\frac{1}{n}\iint\varphi_{d}^2(u)\sum_k\psi_k^2(v)\fx(u)f(u,v)dudv\nonumber\\
&\leq &\frac{\phi_2D_{m_2}}{n}\int\varphi_{d}^2(u)\fx(u)\left(\int f(u,v)dv\right)du\\
&\leq&\frac{\phi_2D_{m_2}}{n}\|\fx\|_\infty.
 \end{eqnarray*}
Hence,
  \begin{eqnarray*}
  \E\sup_{t\in {\mathcal A}}\nu_n^2(t)&\leq&(r+1) \|\fx\|_\infty\frac{\phi_2 D_{m_2}}n
 \end{eqnarray*}
 so that we can take
 $H^2=(r+1)\|\fx\|_\infty\phi_2 D_{m_2}/n$.

 $\bullet$ Thirdly
  \begin{eqnarray*}
{\rm Var}(t(X_1,Y_1))&\leq &\E|t(X_1,Y_1)|^2\\
 &\leq& \iint t^2(u,v)\fx(u) f(u,v)dudv\\
 &\leq & \|t\|_2^2\|f\|_{\infty}\|\fx\|_\infty
  \end{eqnarray*}
 and then  we can take $v=\|f\|_\infty\|\fx\|_\infty $.

Finally
 \begin{align*}
 &\frac{\zeta^2nH^2}{6v}=\frac{\zeta^2(r+1)\phi_2}{6\|f\|_\infty}D_{m_2} \\
 &\frac{\min(\zeta,1)\zeta nH}{21 M}= \frac{\min(\zeta,1)\zeta\sqrt{(r+1)\|\fx\|_\infty}}{21\sqrt{\phi_1}}\sqrt{\frac n{D_{m_1}}}.
 \end{align*}
 According to condition $(CM)$,  on $\Lambda_\rho$,  since $\delta\leq 1$, ${D_{m_1}}\leq \rho n/(\log n)^3$ and $D_{m_2}\geq(\log n)^{2}$.
Thus Talagrand's Inequality gives
\begin{equation*}
 \P_\rho\Big[\sup_{t\in   {\mathcal A}}|\nu_n(t)|\geq (1+2\zeta)H\Big]
 \leq 2\exp(-C\log^{3/2}(n))
  \end{equation*}
with $C$ only depending on $\eta,\zeta,r,\phi_1,\phi_2,\|f\|_\infty,\|\fx\|_\infty$.
Moreover,
  $$(1+2\zeta)H= (1+2\zeta)\sqrt{(r+1)\|\fx\|_\infty\phi_2}\sqrt{\frac{D_{m_2}}{ n}}$$
and, since $\delta >\rho^{-1}\hat{\delta} $ and
$\widehat{\|\fx\|_\infty} >\rho^{-2} \|\fx\|_\infty$
on $\Lambda_\rho$,
  $$ \frac{\sigma(m)}{2\rho^2{\delta}^{-1} \sqrt{\phi_1D_{m_1}}}\geq
 \rho^{-4}(1+\eta)\sqrt{(r+1)\|\fx\|_\infty\phi_2}\sqrt{\frac{D_{m_2}}{n}}.$$
  Then, since $\rho^5=1+\eta$, choosing $\zeta$ such that $1+2\zeta=\rho$ gives
  $$\frac{\sigma(m)}{2\rho^2{\delta}^{-1} \sqrt{\phi_1D_{m_1}}}\geq
  (1+2\zeta)H$$
  and then
  \begin{equation*}
  P_{1,m} \leq
 2\exp(-C\log^{3/2}(n)).
  \end{equation*}

\bigskip

\noindent $[2]$ \emph{Study of $P_{2,m}$}: We now have to bound (with large probability) the term
\begin{eqnarray*}
\|\hat{G}-G\|_{F}^{2}
&= &\sum_{d,d'}\left|\frac{1}{n}\sum_{i=1}^n \varphi_d\varphi_{d'}(X_i)-\E[\varphi_d\varphi_{d'}(X_i)]\right|^2.
\end{eqnarray*}
We use Bernstein's Inequality (Lemma~\ref{bernstein}):
Since $\sup_{u\in\R}|\varphi_d(u)\varphi_{d'}(u)|_\infty\leq \phi_1D_{m_1}$ and
$$\E|\varphi_d\varphi_{d'}(X_1)|^2\leq \iint \varphi_{d}^2\varphi_{d'}^2(u)\fx(u)du\leq \phi_1\|\fx\|_\infty D_{m_1},$$
the assumptions of Lemma~\ref{bernstein} are satisfied with $c=\phi_1D_{m_1}$ and $v=\phi_1\|\fx\|_\infty D_{m_1}$.
If we set $\varepsilon=C_4\sqrt{\frac{D_{m_1}D_{m_2}}{n}}$, with $C_4=(\rho-1)(1+\eta)\sqrt{\phi_1\phi_2\|\fx\|_\infty}/(\rho^3 C_3\delta\sqrt{r+1})$
then, on $\Lambda_\rho$,
$$\varepsilon\leq\frac{(\rho -1)\sigma(m)}{2\rho C_3(r+1)}.$$
Moreover on $\Lambda_\rho$, since $\delta\leq 1$,
 \begin{align*}
 &\frac{n\varepsilon^2}{v}=\frac{C_4^2}{\phi_1\|\fx\|_\infty}D_{m_2}\geq \frac{C_4^2}{\phi_1\|\fx\|_\infty}(\log n)^2\\
 &\frac{n\varepsilon}{c}= \frac{C_4}{\phi_1}\sqrt{\frac{nD_{m_2}}{D_{m_1}}}\geq
 \frac{C_4}{\phi_1\sqrt{\rho\delta}}(\log n)^{5/2}\geq
 \frac{C_4}{\phi_1\sqrt{\rho}}(\log n)^{5/2}.
 \end{align*}
 Then, using Lemma~\ref{bernstein},
\begin{eqnarray*}
\P_\rho\left[\|\hat{G}-G\|_{F}\geq (\rho -1)\frac{\sigma(m)}{2C_3\rho}\right]
 &\leq&\sum_{d,d'}\P_\rho\left(
 \left|\frac{1}{n}\sum_{i=1}^n \varphi_d\varphi_{d'}(X_i)-\E[\varphi_d\varphi_{d'}(X_i)]\right|
 \geq \frac{(\rho -1)\sigma(m)}{2\rho C_3(r+1)}\right)\\
&\leq &  2(r+1)^2 \exp(-C_5\log^{2}(n)),
\end{eqnarray*}
with $C_5$ only depending on $\eta,r,\phi_1,\phi_2,\|f\|_\infty,\|\fx\|_\infty$ and $\delta$. Finally, we denote $$\overline{\mathcal{M}_n}=
\{(m_1,m_2), \quad 2^{m_1}\leq \rho \delta n,\quad D_{m_2}\leq n\}$$
which verifies $\mathcal{M}_n\subset \overline{\mathcal{M}_n}$ on $\Lambda_\rho$.
Gathering all the terms together, we obtain
\begin{eqnarray*}
\P_\rho\left(\exists m\in \mathcal{M}_n\quad \|\hat f_{m}- \mathcal{K}_m(f)\|_{x,2}\geq \frac{\sigma(m)}{2}\right)&\leq&
\P_\rho\left(\exists m\in \overline{\mathcal{M}_n}\quad
\|\hat f_{m}- \mathcal{K}_m(f)\|_{x,2}\geq \frac{\sigma(m)}{2}\right)\\
&\leq&
\sum_{m\in \overline{\mathcal{M}_n}} P_{1,m}+P_{2,m}\\
&\leq& \sum_{m\in \overline{\mathcal{M}_n}}4(r+1)^2 \exp(-C_6\log^{3/2}(n))\\
&\leq & 4(r+1)^2\rho\delta n^2 \exp(-C_6\log^{3/2}(n)),
\end{eqnarray*}
with $C_6$ depending on  $\eta,r,\phi_1,\phi_2,\|f\|_\infty,\|\fx\|_\infty$ and $\delta$, which yields Proposition~\ref{concentrationbis }.

\subsubsection{Proof of Proposition~\ref{omegac}}\label{preuveomegac}
In this Section, we denote
$$\|t\|_X^2:=\int t^2(u)\fx(u)du.$$
We recall that $l_{m_1}$ is the index such that $x$ belongs to the interval $I_{l_{m_1}}$ and as in Section~\ref{preuveconcentration}, we set:
$$F_{m_1}^{(l_{m_1})}={\rm Span}(\varphi_{l_{m_1},d}^m)_{0\leq d\leq r}.$$
We want to bound
\begin{eqnarray*}
 \P(\Omega^c\cap\Lambda_\rho)=\P\left(\exists m_1,\quad  \exists t\in {\rm Span}(\varphi_{l_{m_1},d}^m)_{0\leq d\leq r} \quad\|t\|_n^2< \rho^{-1}\|t\|_X^2 \mbox{ and } \Lambda_\rho\right).
\end{eqnarray*}
Under $(CM)$, we have:
$k_n(r+1)\leq D_{m_1}\leq \hat \delta n /(\log n)^3$, and on $\Lambda_\rho$, we have:
$2^{m_1}\leq \rho \delta n $.
Let $\mu_n$ be the empirical process defined by $$\mu_n(t)=\frac1n\sum_{i=1}^n t(X_i)-\E(t(X_i)).$$ Then, $\mu_n(t^2)=\|t\|_n^2-\|t\|_X^2$, which implies that
$$\P(\Omega^c\cap\Lambda_\rho)\leq \sum_{m_1,2^{m_1}\leq \rho \delta n}\P\left(\sup_{t\in F_{m_1}^{(l_{m_1})}, \|t\|_X=1}|\mu_n(t^2)|>1-\rho^{-1}, \right).$$
But, for all $t\in F_{m_1}^{(l_{m_1})}$ such that $\|t\|_X=1$
$$|\mu_n(t^2)|^2\leq \delta^{-2} \sum_{d,d'} \mu_n^2(\varphi_{l_{m_1},d}^m\varphi_{l_{m_1},d'}^m).$$
Using Lemma~\ref{bernstein}, we easily prove as in Section~\ref{preuveconcentration} that $\forall m\in \mathcal{M}_n$
$$\P\left(| \mu_n(\varphi_{l_{m_1},d}^m\varphi_{l_{m_1},d'}^m)|>(1-\rho^{-1})\delta/(r+1)\right)
\leq 2\exp(-K (\log n)^3)$$
with $K$ depending on $\rho, \phi_1, r, \|\fx\|_\infty,\delta.$
Then
\begin{eqnarray*}
 \P(\Omega^c\cap\Lambda_\rho)&\leq& \sum_{m_1,2^{m_1}\leq \rho \delta n}\P\left( \sum_{d,d'} \mu_n^2(\varphi_{l_{m_1},d}^m\varphi_{l_{m_1},d'}^m)>(\delta(1-\rho^{-1}))^2
 \right)\\
&\leq& \sum_{m_1,2^{m_1}\leq \rho \delta n}\sum_{d,d'}\P\left(  |\mu_n(\varphi_{l_{m_1},d}^m\varphi_{l_{m_1},d'}^m)|>\delta(1-\rho^{-1})/(r+1)\right)\\
&\leq & 2(r+1)^2\sum_{m_1,2^{m_1}\leq \rho \delta n} \exp(-K (\log n)^3)
\leq 2(r+1)^2\rho\delta n\exp(-K \log^3(n)),
\end{eqnarray*}
which yields the result.
\subsubsection{Proof of Proposition~\ref{coarsebound}}\label{preuvecoarsebound}
First, as already noticed, $\|f\|_{x,2}^{2}\leq \|f\|_\infty.$
Now let $m$ be a fixed element of $\mathcal{M}_n$.
Then we denote $l= l_{m_1} $the index such that $x$ belongs to the interval $I_{l}$ and moreover we denote
$$\hat{G}= \hat{G}_m^{(l)},\quad \hat{Z}=\hat{Z}_m^{(l)}, \quad \hat{A}=\hat{A}_m^{(l)},\quad\varphi_d=\varphi_{l,d}^m, \quad \psi_k^m=\psi_k.$$
The elements of $\hat A$  are denoted $(\hat a_{d,k})_{d,k}$  instead of $(\hat a_{(l_{m_1},d),k}^m)_{d,k}$.

If $\Sp(\hat{G})\geq(1+\eta)^{-2/5} \hat{\delta}$ (otherwise $\hat A=0$),
$$\|\hat{G}^{-1}\|_{2}=\rho(\hat{G}^{-1})=(\min (\Sp(\hat{G})))^{-1}\leq (1+\eta)^{2/5}\hat{\delta}^{-1}.$$
Therefore, we have:
\begin{eqnarray*}
 \|\hat A\|_{F}^{2}&\leq&\|\hat{G}^{-1}\|_{2}^{2}\|\hat Z\|_{F}^{2}
 \leq (1+\eta)^{4/5}\hat{\delta}^{-2}\sum_{d,k}\left|\frac{1}{n}\sum_{i=1}^n\varphi_{d}(X_i)\psi_k(Y_i)\right|^{2}\\
&\leq &
 (1+\eta)^{4/5}\hat{\delta}^{-2}\sum_{d,k}\phi_1\phi_2D_{m_1}D_{m_2}\\
&\leq &(1+\eta)^{4/5}\hat{\delta}^{-2}(r+1)\phi_1\phi_2D_{m_1}D_{m_2}^{2}.
\end{eqnarray*}
Finally
\begin{eqnarray*}
\|\hat f_m\|_{x,2}^2
&=&\sum_{k\in K_m}\left(\sum_{d=0}^r\hat a_{dk}\varphi_{d}(x)\right)^2
\leq \|\hat A\|_{F}^2{\phi_1D_{m_1}}\\
 &\leq &(1+\eta)^{4/5}\hat{\delta}^{-2}(r+1)\phi_1^{2}\phi_2D_{m_1}^{2}D_{m_2}^{2}.
 \end{eqnarray*}
\subsubsection{Proof of Theorem~\ref{vitesse}}\label{SecPreuvebiais}
We first assume that $d_1=d_2=1$.
We denote $\mathcal{K}_m^{1}$ the projection on $F_{m_1}$ endowed with the scalar product $(g,h)_X=\int g(z)h(z)\fx(z)dz,$
and $\mathcal{K}_m^{2}$ the projection on $H_{m_2}$ endowed with the usual scalar product $(g,h)_{us}=\int g(z)h(z)dz$.
The projection $\mathcal{K}_m(f)$ can be written for any $u$ and any $y$,
$$\mathcal{K}_m(f)(u,y)=\sum_{k\in K_m} ( f^1(u,.) ,\psi_k^m)_{us} \psi_k^m(y)=\mathcal{K}_m^{2}(f^1(u,.))(y)$$
where $f^1(.,y)=\mathcal{K}_m^{1}(f(.,y)).$
Thus we have the factorization
\begin{eqnarray*}
(\mathcal{K}_m(f)-f)(u,.)
&=&\mathcal{K}_m^{2}(f^1(u,.)-f(u,.))+\mathcal{K}_m^{2}(f(u,.))-f(u,.)
\end{eqnarray*}
and applying Pythagora's theorem
\begin{eqnarray*}
 \|\mathcal{K}_m(f)-f\|^2_{x,2}&= &
 \|\mathcal{K}_m^{2}(f^1(x,.)-f(x,.))\|_{2}^2+\|\mathcal{K}_m^{2}(f(x,.))-f(x,.)\|_{2}^2\\
 &\leq & \|f^1(x,.)-f(x,.)\|_{2}^2+\|\mathcal{K}_m^{2}(f(x,.))-f(x,.)\|_{2}^2.
\end{eqnarray*}
Now, we shall use the following result.
Let $\tau$ be a univariate function belonging to the H\"older space  $\mathcal{H}_{1}(\alpha,L)$ on a interval with length $b$. If $S$ is the space of piecewise polynomials of degree bounded by $r>\alpha-1$ based on the regular partition with $2^J$ pieces, then there exists a constant $C(\alpha,b)$ only depending on $\alpha$ and $b$ such that
$$ d_\infty(\tau, S):=\inf_{t\in S} \|t-\tau\|_{\infty}\leq C(\alpha,b)L2^{-J\alpha}$$
(see for example Lemma 12 in \cite{BBM}).
Let $\mathcal{K}$ the orthogonal projection on $S$ endowed with some scalar product.
We denote $$\vvvert\mathcal{K}\vvvert=\sup_{t\in \L_{\infty}\setminus\{0\}}\frac{\|\mathcal{K}(t)\|_{\infty}}{\|t\|_{\infty}}.$$
Then, for all $t\in S$, since $\mathcal{K}(t)=t$,
\begin{eqnarray*}
\|\tau-\mathcal{K}(\tau)\|_{\infty}=\|\tau-t+\mathcal{K}(t-\tau)\|_{\infty}
\leq (1+\vvvert\mathcal{K}\vvvert)\|t-\tau\|_{\infty}.
\end{eqnarray*}
We obtain:
\begin{eqnarray*}
\|\tau-\mathcal{K}(\tau)\|_{\infty}
\leq (1+\vvvert\mathcal{K}\vvvert )\inf_{t\in S} \|t-\tau\|_{\infty}\leq (1+\vvvert\mathcal{K}\vvvert ) C(\alpha,b)L2^{-J\alpha}.
\end{eqnarray*}
It remains to bound $\vvvert\mathcal{K}\vvvert$ in the following cases.\\

$\bullet$ Case 1: $S$ is the space of piecewise polynomials of degree bounded by $r_1$, endowed with $(.,.)_{X}$ ($S=F_{m_1}$, $\mathcal{K}=\mathcal{K}_m^{1}$).
It is sufficient to apply Lemma~\ref{majbias}
to  the function $\tau(u,y)=t(u)\psi_{k}^m(y)$ to obtain $\vvvert\mathcal{K}\vvvert\leq (r_1+1)\|\fx\|_{\infty}\delta^{-1}$.\\

$\bullet$ Case 2: $S$ is the space of piecewise polynomials of degree bounded by $r_2$, endowed with the usual dot product  ($S=H_{m_2}$, $\mathcal{K}=\mathcal{K}_m^{2}$).
Then it is sufficient to apply the previous case with $\fx$ identically equal to 1, to obtain $\vvvert\mathcal{K}\vvvert\leq (r_2+1)$.\\

Finally, we have obtained the following result:
if $\tau$ is a univariate function belonging to the H\"older space $\mathcal{H}_{1}(\alpha,L)$ then
\begin{eqnarray*}\label{lemholder}
\|\tau-\mathcal{K}_m^1(\tau)\|_{\infty}&\leq &C(\alpha,A,r_1,\|\fx\|_\infty/\delta)LD_{m_1}^{-\alpha},\\
\|\tau-\mathcal{K}_m^2(\tau)\|_{\infty}&\leq &C(\alpha,|B|,r_2)LD_{m_2}^{-\alpha}.
\end{eqnarray*}
Now $f(x,.)$ belongs to the H\"older space $\mathcal{H}_{1}(\alpha_{2},L_{2})$ then
$$\|\mathcal{K}_m^{2}(f(x,.))-f(x,.)\|_{\infty}\leq C_{2}L_{2}D_{m_2}^{-\alpha_2}$$
with $C_{2}$ depending on $\alpha_2,|B|$ and $r_2$.
Moreover, for all $y\in B$, $f(.,y)$ belongs to the H\"older space $\mathcal{H}_{1}(\alpha_{1},L_{1})$ then
$$|f^1(x,y)-f(x,y)|\leq\|\mathcal{K}_m^{1}(f(.,y))-f(.,y)\|_{\infty}\leq C_{1}(\alpha_1,A,r,\|\fx\|_\infty/\delta)L_{1}D_{m_1}^{-\alpha_1}$$
with $C_1$ not depending on $y$.
Finally, since the support of $f(x,.), f^{1}(x,.), \mathcal{K}_m^{2}(f(x,.))$  is compact, we obtain
\begin{eqnarray*}
 \|\mathcal{K}_m(f)-f\|_{x,2}
 &\leq & C_{0}(L_{1}D_{m_1}^{-\alpha_1}+L_{2}D_{m_2}^{-\alpha_2}).
\end{eqnarray*}
with $C_{0}$ depending on $A,|B|,\boldsymbol r,\alpha_1,\alpha_2$ and $\|\fx\|_\infty$ and $\delta$.
We can easily generalize this result to the case $d_1,d_2\geq 2$ and we obtain:
\begin{eqnarray*}
 \|\mathcal{K}_m(f)-f\|_{x,2}
 &\leq & C\sum_{i=1}^d L_i2^{-\alpha_im_i}
\end{eqnarray*}
for $C$  a constant. To conclude, by using Theorem \ref{oracleproba},
it remains to find $(m_1,\ldots,m_d)$ that minimizes
$$(m_1,\ldots,m_d)\longmapsto \sum_{i=1}^d L_i2^{-\alpha_i m_i}+\sqrt{\frac{\prod_{i=1}^d 2^{m_i}}{\delta n}}.$$
Solving this minimization problem shows that $2^{m_i}$ has to be equal to  $L_i^{\frac{1}{\alpha_i}}\Delta_n^{\frac{1}{\alpha_i}}$ up to a constant and
$$\Delta_n=\left(\prod_{i=1}^dL_i^{\frac{1}{\alpha_i}}\right)^{-\frac{\overline{\alpha}}{2\overline{\alpha}+1}}\left(\delta n\right)^{\frac{\overline{\alpha}}{2\overline{\alpha}+1}}.$$
It gives the result.
\appendix
\section{Proofs of technical results}\label{preuvestechniques}
\subsection{Proof of Lemma \ref{tala2}}
We apply the Talagrand concentration inequality given in \citet{kleinrio} to the functions
$s^i(x)=\tau_a(x)-{\mathbb E}(\tau_a(U_i))$
and we obtain
$$\P(\sup_{a\in A}|\nu_n(a)|\geq H+\lambda)\leq 2\exp\left(-\frac{n\lambda^2}{2(v+4HM)+6M\lambda}\right).$$
Then we modify this inequality following \citet{BM98} Corollary 2 p.354. It gives
\begin{equation}
  \label{talapr}
  \P(\sup_{a\in\mathcal A}|\nu_n(a)|\geq (1+\zeta)H+\lambda)\leq 2\exp\left(-\frac{n}{3}\min
\left(\frac{\lambda^2}{2v},\frac{\min(\zeta,1)\lambda}{7M}\right)\right).
\end{equation}
To conclude, we set $\lambda=\zeta H$.
\subsection{Proof of Lemma \ref{lemmafx}}
The lemma is a consequence of (\ref{Cfx}) used with $\lambda=\rho-1$, $\lambda=1-\rho^{-1}$, $\lambda=\rho^2-1$ or $\lambda=1-\rho^{-2}$.
Indeed, under $(\ref{Cfx})$, with probability $1-\kappa\exp(-(\log n)^{3/2})$, for all $t\in V_n(x)$,
$|\fx(t)-\hfx(t)|\leq \lambda  |\hfx(t)|$, which implies
$$(1-\lambda)|\hfx(t)|\leq |\fx(t)|\leq (1+\lambda ) |\hfx(t)|$$
and then
$$(1+\lambda)^{-1}|\fx(t)|\leq |\hfx(t)|\leq (1-\lambda)^{-1} |\fx(t)|.$$
Thus, with probability $1-\kappa\exp(-(\log n)^{3/2})$,  $(1+\lambda)^{-1}\delta\leq \hat\delta\leq (1-\lambda)^{-1} \delta$
and $(1+\lambda)^{-1}\|\fx\|_\infty\leq \|\hfx\|_\infty\leq (1-\lambda)^{-1} \|\fx\|_\infty.$
\subsection{Proof of Lemma \ref{sup1}}
We have:
\begin{eqnarray*}
\|f_1*f_2\|^2_{x,2}&=&\int (f_1*f_2)^2(x,y)dy
=\int\left(\iint f_1(x-u,y-v)f_2(u,v)dudv\right)^2dy\\
&\leq&\int\left(\iint f_1^2(x-u,y-v)|f_2(u,v)|dudv\times\iint |f_2(u,v)|dudv\right)dy\\
&=&\|f_2\|_1\iint\|f_1(x-u,\cdot)\|_2^2|f_2(u,v)|dudv
\leq\sup_{t\in V_n(x)} \|f_1\|_{t,2}^2\times \|f_2\|_1^2.
\end{eqnarray*}
\subsection{Proof of Lemma \ref{majbias}}
Let $l$ the index such that $x$ belongs to the interval $I_l$.
We denote
$$\varphi_d=\varphi_{l,d}^m, \quad \psi_k=\psi_k^{m},$$
$$I(\tau)=\left(\iint
\varphi_{d}(u)\psi_k(y)\tau(u,y)\fx(u)dudy\right) _{0\leq d\leq r, k \in K_m}$$
and
$$\mathcal{K}_m(\tau)(x,y)=\sum_{k}\sum_{d}b_{d,k}\varphi_{d}(x)\psi_{k}(y).$$
Lemma~\ref{obvious} shows that the matrix of coefficients  $B=(b_{d,k}) _{0\leq d\leq r, k \in K_m}$  verifies the equation $GB=I(\tau),$
with
$${G}=\E(\hat{G})=\left(\langle
\varphi_{d_1},\varphi_{d_2}
\rangle_{X}\right)
_{0\leq d_1, d_2\leq r}.$$
Now, using \eqref{min1},
\begin{eqnarray*}
\| \mathcal{K}_m(\tau)\|_{x,2}^{2}&=&\sum_{k}(\sum_{d}b_{d,k}\varphi_{d}(x))^{2}
\leq \sum_{d}\varphi_{d}^{2} (x) \sum_{k}\sum_{d}b_{d,k}^{2}\\
&\leq & \phi_1 D_{m_1}\|B\|^{2}_{F}\leq  \phi_1 D_{m_1}\|G^{-1}\|_{2}^{2}\|I(\tau)\|^{2}_{F}\\
&\leq &  \phi_1 D_{m_1}\delta^{-2}\|I(\tau)\|^{2}_{F}.
\end{eqnarray*}
Now we denote $\text{Proj}_{H_{m_2}}$ the usual $\L_2$ orthogonal projection on $H_{m_2}$ and $(\cdot,\cdot)_{us}$  the standard $\L_2$ dot product.
Notice that for any function $s\in \L_2(\R)$, $\sum_{k\in K_m}( s, \psi_k)_{us}^2=\int |\text{Proj}_{H_{m_2}}(s)|^{2}(y)dy\leq\int s^{2}(y)dy$. Then
\begin{eqnarray*}
\|I(\tau)\|_{F}^{2}
&=& \sum_{d=0}^r\sum_{k\in K_m}( \int \varphi_{d}(u)\tau(u,.)\fx(u)du , \psi_k)_{us}^2\\
&\leq &
\sum_{d=0}^r\int \left(\int \varphi_{d}(u)\tau(u,y)\fx(u)du\right)^{2}dy
 \leq\int \left(\int_{I_l} \tau^2(u,y)\fx^2(u)du\right)dy
 \end{eqnarray*}
 using that for any function $s$,
$\sum_d (\int s \varphi_d)^2\leq\int_{I_l} s^2$. Next,
using that $I_{l}$ is an interval with length $4A(r+1)D_{m_1}^{-1}$,
\begin{eqnarray*}
\|I(\tau)\|_{F}^{2}
&\leq &\sup_{t\in I_{l}}\|\tau\|_{t,2}^{2}\int_{I_l}\fx^2(u)du
\leq
4A(r+1)D_{m_1}^{-1}\|\fx\|_{\infty}^{2}\sup_{t\in I_{l}}\|\tau\|_{t,2}^{2}.
\end{eqnarray*}
Finally
\begin{eqnarray*}
\| \mathcal{K}_m(\tau)\|_{x,2}^{2}
&\leq &  \phi_1 D_{m_1}\delta^{-2}4A(r+1)D_{m_1}^{-1}\|\fx\|_{\infty}^{2}\sup_{t\in I_{l}}\|\tau\|_{t,2}^{2}\\
&\leq & (r+1)^{2}\|\fx\|_{\infty}^{2} \delta^{-2}\sup_{t\in V_{n}(x)}\|\tau\|_{t,2}^{2}
\end{eqnarray*}
and the lemma is proved.
\section{Tables for simulation results}\label{simustables}
In this appendix, for each example and each procedure, we give the approximated mean squared error based on $N=100$ samples for different values of $n$,  different values of the parameter $\eta$ and different values of $x$. We give in bold red the minimal value of the approximated mean squared error. For the kernel estimator and Examples 1 and 2, we distinguish the case where $\fx$ is known or not.
\begin{table}[!htp]\begin{center}
\begin{tabular}{|c|c|c|c|c|c|c|c|c|c|c|}\hline
Ex 1&\multicolumn{5}{|c|}{$\fx$ known}& \multicolumn{5}{|c|}{$\fx$ unknown} \\\hline
  $\eta$     & $-0.2$  & $0.5$& $1$ &$2$      & $3$  &  $-0.2$ & $0.5$& $1$ & $2$ & $3$\\\hline
$n=250$      &  1.285  &0.061  &\red{0.017}& 0.020 & 0.029&1.368 &0.033 &\red{0.028} & 0.042 &0.062 \\\hline
$n=500$      &  0.673  &0.019  &\red{0.009}& 0.010 &0.018 &0.685 &0.016 &\red{0.009} & 0.011 &0.018\\\hline
$n=1000$     &  0.336  & 0.013 &\red{0.006}& \red{0.006} &0.009 &0.329 &0.013 &\red{0.006} & 0.007 &0.010\\\hline
\end{tabular}
\caption{Mean squared error for the kernel estimator at $x=0.5$ for Example 1}
\label{table1.simul}
\end{center}
\end{table}
\begin{table}[!htp]\hspace{6.7cm}\begin{center}
\begin{tabular}{|c|c|c|c|c|c|}\hline
Ex 1& \multicolumn{5}{|c|}{$\fx$ unknown} \\
\hline
$\eta$ & $-0.2$ &  $0.5$& $1$ & $2$ & $3$\\\hline
$n=250$&    0.492  &  \red{0.192} &  0.222 &  0.232 & 0.231\\\hline
$n=500$&    0.087  &  \red{0.076}  &  0.119 &  0.211 & 0.229\\\hline
$n=1000$&    0.051  &  \red{0.047}  &  0.055 & 0.070 &  0.138 \\\hline
\end{tabular}
\caption{Mean squared error for the projection estimator at $x=0.5$ for Example 1}
\label{table1}
\end{center}
\end{table}
\begin{table}[!htp]\begin{center}
\begin{tabular}{|c|c|c|c|c|c|c|c|c|c|c|}\hline
Ex 2 &\multicolumn{5}{|c|}{$\fx$ known}& \multicolumn{5}{|c|}{$\fx$ unknown} \\\hline
$\eta $ &  $-0.2$  & $0.5$ & $1$  &2   & $3$ &    $-0.2$ & $0.5$  &$1$    &2   & $3$\\\hline
$n=250$  &  0.038  &0.008  &\red{0.006} & 0.007 &0.009 &    0.042  &  0.008 & \red{0.006} & 0.008   & 0.009 \\\hline
$n=500$ &  0.021   &0.006  &\red{0.004} & 0.005 &0.006 &     0.025 &0.006   &\red{0.004}  & 0.005   &0.007\\\hline
$n=1000$& 0.01     & 0.004 &\red{0.003} & 0.004 &0.005 &    0.012  &0.004   &\red{0.003}  & 0.004   &0.005\\\hline
\end{tabular}
\caption{Mean squared error for the kernel estimator  at $x=0.5$ for Example 2}
\label{table3.simul}
\end{center}
\end{table}
\begin{table}[!htp]\hspace{6.7cm}\begin{center}
\begin{tabular}{|c|c|c|c|c|c|}\hline
Ex 2& \multicolumn{5}{|c|}{ $\fx$ unknown} \\
\hline
$\eta$ & $-0.2$  & $0.5$& $1$ & 2 & $3$\\\hline
$n=250$&   0.154  &  \red{0.104}   & 0.128  & 0.152 &  0.158 \\\hline
$n=500$&   \red{0.064}   & 0.070  &  0.090 &0.103 &   0.123\\\hline
$n=1000$ &  \red{0.047}  &  0.060   & 0.063 &0.074 &  0.088
\\\hline
\end{tabular}
\caption{Mean squared error for the projection estimator at $x=0.5$ for Example 2}
\label{table3}
\end{center}
\end{table}
\begin{table}[!htp]\begin{center}
\begin{tabular}{|c|c|c|c|c|c|c|}\hline
\multicolumn{2}{|c|}{Ex 3 }  & \multicolumn{5}{|c|}{$\fx$ unknown} \\\hline
        &$x \backslash\eta$  &  $-0.2$   & $0.5$ & $1$  &2  & $3$\\\hline\hline
        &0                   &0.514	&0.016&	0.013&\red{0.012}&	0.019\\\cline{2-7}
$n=250$  &0.36               & 	0.092	&\red{0.062}&	0.080&0.112&	0.134 \\\cline{2-7}
        &1                   &1.709&	0.015	&\red{0.009}&\red{0.009}&	0.016  \\\hline\hline
        &0                   & 	0.269&	0.013&	0.013&	\red{0.009}&0.010  \\\cline{2-7}
$n=500$  &0.36                 & 	0.109&	0.040&	\red{0.039}&0.063&	0.094\\\cline{2-7}
        &1                   & 	0.601&	0.010&	0.009&\red{0.006}&	0.008\\\hline\hline
        &0                   &	0.126&	0.011	&0.011&0.008&	\red{0.006} \\\cline{2-7}
$n=1000$ &0.36                & 	0.104&	0.029&	\red{0.024}&0.037&	0.056 \\\cline{2-7}
        &1                   & 	0.265&	0.006&	0.007&	\red{0.004}&\red{0.004}  \\\hline
\end{tabular}
\caption{Mean squared error for the kernel estimator  at $x=0$, $x=0.36$ and $x=1$ for Example 3}
\label{table2.simul}
\end{center}
\end{table}
\begin{table}[!htp]\hspace{6.3cm}\begin{center}
\begin{tabular}{|c|c|c|c|c|c|c|}\hline
\multicolumn{2}{|c|}{Ex 3 } & \multicolumn{5}{|c|}{$\fx$ unknown} \\
\hline
        & $x\backslash\eta$ & $-0.2$   & $0.5$  & $1$    &2 & $3$\\\hline
        &0                  &   \red{0.029}  & 0.035 &  0.041 & 0.051   & 0.060\\\cline{2-7}
$n=250$  &0.36               &   0.186  &  0.188 & 0.183 & 0.172   & \red{0.170} \\\cline{2-7}
        &1                  &   \red{0.033}  &  0.038 &  0.044 & 0.064   & 0.099\\\hline\hline
        &0                  &   \red{0.020}  &  0.028 &  0.033 & 0.036  &   0.038\\\cline{2-7}
$n=500$  &0.36               &  \red{0.169}  &  0.184 &  0.177  & 0.172  &  0.170\\\cline{2-7}
        &1                  &   \red{0.027}  &  0.029 &  0.030  & 0.032 &  0.035\\\hline\hline
        &0                  &   \red{0.012}  & 0.018 &  0.023 & 0.031   & 0.034\\\cline{2-7}
$n=1000$ &0.36               & \red{0.160} &  0.161 &   0.166  & 0.170    & 0.169\\\cline{2-7}
        &1                  & \red{0.023}   & 0.025 &  0.028  & 0.029   & 0.028\\\hline
\end{tabular}
\caption{Mean squared error for the projection estimator at $x=0$, $x=0.36$ and $x=1$ for Example~3}
\label{table2}
\end{center}
\end{table}
\begin{table}[h]\begin{center}
\begin{tabular}{|c|c|c|c|c|c|c|}\hline
\multicolumn{2}{|c|}{Ex 4 }  & \multicolumn{5}{|c|}{$\fx$ unknown} \\\hline
        &$x \backslash\eta$   & $-0.2$   & $0.5$ & $1$  &2  & $3$\\\hline
&0&		0.016&	\red{0.007}	&\red{0.007}&	0.009&0.013\\\cline{2-7}
n=250&0.36&	0.082	&\red{0.03}	&0.037&	0.048&0.055\\\cline{2-7}
&1&	 	0.026	&\red{0.006}&	\red{0.006}&	0.009&0.0119\\\hline\hline
&0&		0.009&	\red{0.004}&	\red{0.004}&	0.006&0.009\\\cline{2-7}
n=500&0.36&		0.057&	\red{0.019}&	0.023&0.034	&0.043\\\cline{2-7}
&1&		0.016&	\red{0.005}&	\red{0.005}&0.006	&0.008\\\hline\hline
&0&		0.004&	\red{0.003}&	\red{0.003}&	0.004&0.005\\\cline{2-7}
n=1000&0.36&	 0.037&	\red{0.013}&	0.014&	0.021&0.03\\\cline{2-7}
&1 &	0.008&	\red{0.003}&	\red{0.003}&	0.004&0.005\\\hline
\end{tabular}
\caption{Mean squared error for the kernel estimator  at $x=0$, $x=0.36$ and $x=1$ for Example 4}
\label{table4.simul}
\end{center}
\end{table}
\begin{table}[h]\hspace{6.3cm}\begin{center}
\begin{tabular}{|c|c|c|c|c|c|c|}\hline
\multicolumn{2}{|c|}{ Ex 4} & \multicolumn{5}{|c|}{$\fx$ unknown} \\
\hline
        & $x\backslash\eta$ &  $-0.2$  & $0.5$    & $1$   &2 & $3$\\\hline
        &0                  &    \red{0.028} &   0.030 & 0.032 & 0.036 & 0.040\\\cline{2-7}
$n=250$  &0.36               &   0.103  & 0.102 &  0.099  &0.096&  \red{0.095}\\\cline{2-7}
        &1                  &    \red{0.030} &   0.036 &  0.038  & 0.049 & 0.066 \\\hline\hline
        &0                  &   \red{0.022}  &  0.024 &  0.024 & 0.029&  0.032 \\\cline{2-7}
$n=500$  &0.36               &   0.098  &  0.099 &  0.097& 0.094 &  \red{0.094} \\\cline{2-7}
        &1                  &  \red{0.026}   & 0.027 &  0.028  &0.033 & 0.036\\\hline\hline
        &  0                &    \red{0.020} &   0.020 &  0.021 & 0.021 & 0.023\\\cline{2-7}
$n=1000$ &0.36               &  \red{0.082}   & 0.083 &  0.093  & 0.095 & 0.094\\\cline{2-7}
        &1                  &   0.023  &  0.023   & \red{ 0.022} & 0.026 &   0.028\\\hline
\end{tabular}
\caption{Mean squared error for the projection estimation in $x=0$, $x=0.36$, $x=1$ for Example 4}
\label{table4}
\end{center}
\end{table}

\newpage	

\noindent{\bf Acknowledgements:} The research of Claire Lacour and Vincent Rivoirard is partly supported by the french Agence Nationale de la Recherche (ANR 2011 BS01 010 01 projet Calibration). Karine Bertin has been partially supported by Project ECOS-CONICYT C10E03 and by the grant ANILLO ACT--1112, CONICYT-PIA, Chile. The authors wish to thank two anonymous referees who each made helpful suggestions that improved the presentation of the paper.
\bibliographystyle{apalike}
\bibliography{biblio}

\end{document}